\documentclass[preprint,12pt]{elsarticle}

\usepackage{amssymb}
\usepackage{amsmath}

\usepackage{bm}
\usepackage{url}
\usepackage{algorithm}
\usepackage{algorithmic}
\usepackage{caption}
\usepackage{subcaption}
\usepackage{tabularx}
\usepackage{multirow}
\usepackage{geometry}
\usepackage{pdflscape}
\usepackage{hyperref}

\usepackage{amsthm}
\newtheorem{remark}{Remark}

\journal{Aerospace Science and Technology}

\begin{document}

\begin{frontmatter}



\title{Eigendecomposition Parameterization of Penalty Matrices for Enhanced Control Design: Aerospace Applications}




\author[aff1]{Nicholas P. Nurre}
\author[aff1]{Ehsan Taheri\corref{cor1}}
\ead{ezt0028@auburn.edu}

\cortext[cor1]{Corresponding author}

\address[aff1]{Department of Aerospace Engineering, Auburn University, 141 Engineering Dr, Auburn, AL 36849, USA}

\begin{abstract}
Modern control algorithms require tuning of square weight/penalty matrices appearing in quadratic functions/costs to improve performance and/or stability output. Due to simplicity in gain-tuning and enforcing positive-definiteness, diagonal penalty matrices are used extensively in control methods such as linear quadratic regulator (LQR), model predictive control, and Lyapunov-based control. In this paper, we propose an eigendecomposition approach to parameterize penalty matrices, allowing positive-definiteness with non-zero off-diagonal entries to be implicitly satisfied, which not only offers notable computational and implementation advantages, but broadens the class of achievable controls. We solve three control problems: 1) a variation of Zermelo’s navigation problem, 2) minimum-energy spacecraft attitude control using both LQR and Lyapunov-based methods, and 3) minimum-fuel and minimum-time Lyapunov-based low-thrust trajectory design. Particle swarm optimization is used to optimize the decision variables, which will parameterize  the penalty matrices. The results demonstrate improvements of up to 65\% in the performance objective in the example problems utilizing the proposed method.
\end{abstract}



\begin{keyword}
Gain tuning \sep Penalty matrix parameterization \sep Eigendecomposition \sep LQR \sep Attitude control \sep Low-thrust \sep Spacecraft \sep Trajectory \sep Lyapunov \sep Quadratic cost \sep MPC
\end{keyword}

\end{frontmatter}



\section{Introduction}
Modern control algorithms rely on tuning a set of gains to achieve a desired stability and/or performance output. For instance, positive-definite and/or positive semi-definite penalty matrices are used to construct quadratic costs or candidate quadratic control Lyapunov functions. Diagonal parameterization of the penalty (also referred to as the ``weighting'' or ``parameter'') matrices is used extensively in LQR controllers \cite{zhou_non-linear_2005,marco_automatic_2016}, model-predictive control (MPC) \cite{dicairano_model_2012,nguyen_tinympc_2024}, sliding mode control \cite{crassidis_sliding_1996}, control Lyapunov functions (CLFs) \cite{petropoulos_low-thrust_2004}, and control barrier function methods \cite{ames_control_2017}. In this work, we consider quadratic CLFs and LQR controllers and investigate the impact of various parameterizations of the penalty matrices on control performance. Lyapunov control (LC) laws derived from CLFs can generate near-optimal control solutions and are used in many aerospace applications. In astrodynamics and spacecraft trajectory optimization, some well-known low-thrust LC laws include Q-law \cite{petropoulos_low-thrust_2004} and the Chang-Chichka-Marsden (CCM) control law \cite{eui_chang_lyapunov-based_2002}. LC laws are also used for spacecraft attitude control \cite{junkins_introduction_2000,schaub_analytical_2018}. LQR is used extensively for solving spacecraft attitude control problems \cite{zhou_non-linear_2005}.

Diagonal parameterization of the penalty matrices is quite common due to 1) implementation simplicity, 2) ease of enforcement of sign-definiteness by constraining the signs of the diagonal entries, and 3) making the gain-tuning intuitive since diagonal entries correspond to each state and control without any cross-coupling. However, diagonal matrices do not span the full solution space, as they neglect cross-coupled terms, leading to suboptimal gains (in the sense of some particular performance and/or stability criterion). Thus, it will be advantageous to consider penalty matrices with non-zero off-diagonal terms, which we refer to herein as ``full'' matrices. This raises a key question: How can a sign-definite full matrix be efficiently parameterized? A straightforward approach is eigendecomposition \cite{strang_introduction_2023}, as proposed in Ref. \cite{junkins_introduction_2000}. Sign-definiteness is trivially enforced through the signs of the eigenvalues. However, parameterizing the associated eigenvector rotation matrices remains less straightforward.

Multiple methods exist for parameterizing orthogonal matrices. Ref. \cite{shepard_representation_2015} outlines four parameterizations that are generalizable to \(N\times N\) dimensions and require the minimum number of design variables of \(M=N(N-1)/2\). One of the first parameterizations discovered, the Cayley transform \cite{cayley_motion_1843,bar-itzhack_minimal_1990,schaub_principal_1995,junkins_introduction_2000,shepard_representation_2015,schaub_analytical_2018}, has been found useful in various practical engineering applications such as in attitude kinematic representations \cite{schaub_principal_1995,schaub_analytical_2018} and efficient solution of the matrix Riccati differential equation \cite{oshman_eigenfactor_1985,bar-itzhack_minimal_1990}. Ref. \cite{nurre_expanding_2024} proposes another parameterization of orthogonal matrices similar to the work of Ref. \cite{hoffman_generalization_1972} requiring \(M\) decision variables and is based on successively projecting unit vectors (starting with an \(N\)-dimensional unit vector constructed with spherical coordinates \cite{blumenson_derivation_1960}). Our goal is to parameterize square penalty matrices as efficiently as possible over the entire solution space, by decomposing the matrices through their 1) eigenvalues and 2) eigenvector matrices. Choosing how to parameterize the orthogonal eigenvector matrix can have a substantial impact on the penalty matrices and the search for more optimal control gains. While full penalty matrices can improve controller performance, a drawback is knowing how to choose the increased number of free variables (\(N\) for diagonal matrices versus \((N^2+N)/2\) for full matrices) which also may not necessarily be physically meaningful. 

Traditional gain-tuning techniques include the Ziegler--Nichols step response method \cite{ziegler_optimum_1942,astrom_revisiting_2004} for PID gains and Bryson's rule (stemming from Section 5.4 Example 2 of Ref. \cite{bryson_applied_1975}) which serves as a starting point for manual trial-and-error gain tuning of LQR controllers. However, modern stochastic/evolutionary optimization techniques and computing technology have led to improved gain-tuning methods. For example, Ref.~\cite{trimpe_self-tuning_2014} uses a stochastic optimization algorithm that approximates the gradient of the performance metric to optimize LQR gains based on a performance metric evaluated from experiment. The same gains are improved further in Ref. \cite{marco_automatic_2016} by using global Bayesian optimization which better utilizes the experimental data. Ref.~\cite{lee_design_2005} uses genetic algorithm and simulated annealing to find the best diagonal penalty matrix for Q-law to form a time-of-flight and propellant mass Pareto front for the same low-thrust transfer cases we solved in Ref.~\cite{nurre_expanding_2024} and in this paper. Automated gain tuning will be necessary for optimizing the larger number of decision variables used to parameterize the full penalty matrices. We elect to use the metaheuristic optimization method particle swarm optimization (PSO) \cite{kennedy_particle_1995} to tune the gains in this work for its simplicity and because it is already available in the MATLAB Global Optimization Toolbox.

The main contribution of this work is to extend the work in Ref. \cite{nurre_expanding_2024} by considering additional examples and analyses for the purpose of investigating the improvement in performance of feedback control laws with full penalty matrices in a variety of control methodologies and aerospace applications. We parameterize the full penalty matrices with eigendecomposition, which allows for the straightforward enforcement of positive-definiteness. The Cayley transform \cite{shepard_representation_2015}, Givens rotation \cite{shepard_representation_2015}, and a parameterization based on generalized Euler angles and Gram-Schmidt process (GEAGSP) \cite{nurre_expanding_2024} are presented, but only the GEAGSP parametrization is used to obtain the results because of the existing code architecture. There may be certain benefits in using one orthogonal matrix parameterization over another due to factors such as computational efficiency \cite{shepard_representation_2015}. The decision variables of the parameterization methods are optimized with the metaheuristic particle swarm optimization (PSO) algorithm \cite{kennedy_particle_1995}. This choice of optimization method was arbitrary as the optimization of the decision variables can be achieved by any meta-heuristic optimization method and the optimization algorithm is not the focus of this study. As such, parameters of the PSO algorithm such as swarm size were also arbitrarily chosen to achieve (what we qualitatively perceive to be) optimal solutions. In the present work, we focus on investigating the advantages of using full penalty matrices for solving three control problems: 1) a variation of Zermelo’s navigation problem, 2) minimum-energy spacecraft attitude control using both LQR and Lyapunov-based methods, and 3) minimum-fuel and minimum-time Lyapunov-based low-thrust trajectory design. To our knowledge, no work explores using a full penalty matrix in aerospace control applications for the sake of improving a certain performance objective such as propellant consumption, aside from our earlier papers \cite{nurre_expanding_2024,nurre2024EclipseFeasible} (Ref. \cite{junkins_simultaneous_1985} considers full parameter matrices for the purpose of designing a flexible structure and its respective controller simultaneously based on eigenvalue placement). Consideration of the off-diagonal penalty terms, as it is shown herein, can enlarge the solution set and result in improved controller performance and reduced values of the cost function for each example problem. Our hope is that exemplifying and quantifying this improvement for a variety of aerospace control problems along with providing a way to efficiently parameterize full penalty matrices will potentially allow practitioners to explore improving their own control algorithms.

The paper is organized with Section \ref{sec: parametrization} presenting a review of the parameterization methods of positive-definite matrices. Sections \ref{sec: zermelo}, \ref{sec: spacecraft attitude control}, and \ref{sec: low thrust trajectory optimization} present the comparisons between the two types of penalty matrices for a variation of Zermelo’s navigation problem, minimum-energy spacecraft attitude control using both LQR and Lyapunov-based methods, and minimum-fuel and minimum-time Lyapunov-based low-thrust trajectory design, respectively. Section \ref{sec: conclusion} concludes the paper.

\section{Penalty Matrix Parameterization} \label{sec: parametrization}

One approach to parameterize positive-definite penalty matrices is to consider a  symmetric matrix and thus use \(M=(N^2+N)/2\) decision variables for parameterization. Let \(\bm{K}\) denote a generic penalty matrix and let \(\bm{k}=(k_1,\dots,k_N,\dots,k_M)\), the matrix can be written as,
\begin{equation}
    \bm{K} = \begin{bmatrix}
    k_1  & k_2     & k_3  & \cdots & k_N \\
    k_2  & k_{N-1} & k_{N-2}  & \cdots & k_{2N-1} \\
    k_3  & k_{N-2}  & k_{2N-1}    & \cdots & k_{3N-2} \\
    \vdots & \vdots & \vdots & \ddots & \vdots \\
    k_N & k_{2N-1} & k_{3N-1} & \cdots & k_M
    \end{bmatrix},
\end{equation}
where the diagonal parameters, \(\bm{k}_\text{diag}\), and the off-diagonal parameters, \(\bm{k}_\text{off-diag}\), are bounded as \(\bm{k}_\text{diag}\in\mathbb{R}^+\) (leveraging the fact that positive-definite matrices always have positive diagonal elements) and \(\bm{k}_\text{off-diag}\in\mathbb{R}\). Then at each iteration within an optimization process, we can check if \(\bm{K}\) is positive-definite in a number of ways. For example, we can calculate its eigenvalues, but this can become computationally demanding as the number of dimensions increases. While this simple parameterization of \(\bm{K}\prec 0\) only requires the minimum amount of variables needed to parameterize a positive-definite matrix, \(M\), the search process over \(\bm{k}\) can be made more efficient by only considering a search space that will always result in a positive-definite \(\bm{K}\). This could be achieved by explicitly enforcing positive-definiteness through nonlinear constraints on the decision variables, \(\bm{k}\). However, this method also has the limitations of optimization with nonlinear constraints, which won't scale appropriately to higher dimensions.

A preferred approach is to implicitly satisfy the sign-definiteness property of the penalty matrix by using eigendecomposition \cite{junkins_introduction_2000,strang_introduction_2023}. Here, we can parameterize \(\bm{K}\) as,
\begin{equation} \label{eq: eigendecompisition}
    \bm{K} = \bm{Q} \bm{\Lambda} \bm{Q}^{\top},
\end{equation}
where \(\bm{\Lambda}=\text{diag}\left(\lambda_1,\dots,\lambda_N\right)\) are the eigenvalues and \(\bm{Q}\in\text{O}(N)\) denotes the orthogonal matrix of eigenvectors. The eigenvalues and eigenvectors can be thought of as representing the scale/size deformation and the rotational deformation, respectively, of \(\bm{K}\) \cite{junkins_introduction_2000}. Positive-definiteness is implicitly enforced by ensuring all eigenvalues are positive. 

We seek to parameterize \(\bm{Q}\) with the minimum number of parameters needed, \(M=(N^2-N)/2\), in order to span \(\text{O}(N)\), i.e., we desire \(\bm{Q}(\bm{\phi}) \in\text{O}(N)\) where \(\bm{\phi}\) denotes the set of parameters with \(\text{dim}(\bm{\phi})=M\). Because we use derivative-free optimization exclusively in this work, we do not need \(\partial \bm{Q}/\partial \bm{\phi}\). However, gradient-based optimization may further improve results (see Ref. \cite{shepard_representation_2015} for a detailed discussion on the computational cost of gradient calculations for each orthogonal-matrix parameterization presented). 

We outline three (out of the many) possible parameterizations methods: Cayley transform (CT), Givens Rotation (GR), and one based on generalized Euler angles and Gram-Schmidt Process (GEAGSP). We note that only the GEAGSP parameterization was used in all the results in this paper. Through numerical simulations, we have found that the CT and GR parameterizations to be comparable in performance. However, there may be benefits of using one parameterization over another for reasons such as computational efficiency, properties of the inverse mapping if required, or nonlinearity of the search space. Thus, the best parameterization could be explored further in future work.

\subsection{Cayley Transform}

The Cayley transform can be used to parameterize the special orthogonal matrices, \(\text{SO}(N)\) \cite{junkins_introduction_2000,shepard_representation_2015}. Letting \(\bm{X}\) denote an \({N\times N}\) skew-symmetric matrix, then the orthogonal matrix \(\bm{Q}(\bm{X}(\bm{\phi}))\) is defined with the Cayley transform as,
\begin{equation}
    \bm{Q} = \left(\bm{I} + \bm{X}\right)\left(\bm{I} - \bm{X}\right)^{-1}.
\end{equation}
Please refer to Ref. \cite{shepard_representation_2015} for more details such as the reverse transformation and its smoothness and continuity properties.

\subsection{Givens Rotation}
The Givens rotation parameterization from Ref. \cite{shepard_representation_2015} can be thought of as a succession of plane transformations. Letting \(\bm{G}(i,j,\theta)\) denote each plane transformation defined as,
\begin{equation}
    \mathbf{G}(i,j,\theta) =
    \begin{bmatrix}
    1 & \cdots & 0 & \cdots & 0 & \cdots & 0 \\
    \vdots & \ddots & \vdots & & \vdots & & \vdots \\
    0 & \cdots & \cos\theta & \cdots & -\sin\theta & \cdots & 0 \\
    \vdots & & \vdots & \ddots & \vdots & & \vdots \\
    0 & \cdots & \sin\theta & \cdots & \cos\theta & \cdots & 0 \\
    \vdots & & \vdots & & \vdots & \ddots & \vdots \\
    0 & \cdots & 0 & \cdots & 0 & \cdots & 1
    \end{bmatrix},
\end{equation}
where the \(\sin{\theta}\) and \(\cos{\theta}\) populate the elements intersected by the \(i\)-th column and \(j\)-th row. The diagonal elements are all 1 except when a trig function populates the diagonal. More formally, the nonzero elements are defined as,
\begin{align}
    G(i,j,\theta)_{ii} & = G(i,j,\theta)_{jj} = \cos\theta, \\ 
    G(i,j,\theta)_{ji} & = -G(i,j,\theta)_{ij} = \sin\theta, \quad \text{for } i > j, \\
    G(i,j,\theta)_{kk} & = 1 \quad \text{for } k \notin \{i, j\}.
\end{align}
Let $\bm{\theta} = [\theta_1,\cdots,\theta_M]$, \(\bm{Q}(\bm{G}(\bm{\theta}))\) is defined as,
\begin{equation}
    {\bm{Q}} = {\bm{G}}_1^\top \dots {\bm{G}}_M^\top.
\end{equation}
\begin{remark}
    It is well known that these plane rotations are ambiguous at certain angles, since there are multiple sets of angles that can represent the same rotation matrix. In 3 dimensions, this is known as ``gimbal lock.'' We note this introduces redundancies, however, we find this characteristic to be ignorable when using heuristic optimization to tune the parameters. This could also be avoided in the optimization algorithm by detecting when one angle reaches an ambiguous value and setting the rest of the angles to arbitrarily fixed values \cite{shepard_representation_2015}.
\end{remark}

\subsection{Generalized Euler Angles and Gram–Schmidt Process}

Another approach to parameterize the \(\bm{Q}\) matrix is to use a generalization of Euler angles derived in \cite{blumenson_derivation_1960} and perform successive orthogonalization through the Gram-Schmidt process, as we proposed in Ref. \cite{nurre_expanding_2024}. This method is similar to the earlier work of \cite{hoffman_generalization_1972}. Let \(\theta_{ij}~\forall i=1,\dots,N,j=1,\dots,N-1\), denote the angle-like parameters of \(\bm{Q}\), this parameterization can be described by Algorithm \ref{alg: rotation} where 
\begin{equation}    
    \bm{v}_i^p = \begin{bmatrix}
        \cos(\theta_{i1}) \\ \sin(\theta_{i1}) \cos(\theta_{i2}) 
        \\ \vdots \\ \prod_{j=1}^{N-i-1} \sin(\theta_{ij}) \cos{(\theta_{i(N-i)})}\\
        \prod_{j=1}^{N-i} \sin(\theta_{ij})
    \end{bmatrix},
\label{eq:unit_vec}
\end{equation}
where \(\bm{v}_i^p\) is a unit vector in a \((N-i + 1)\)-dimensional space. The superscript \(p\) denotes that \(\bm{v}_i\) is not expressed in the original \(N\)-dimensional space and it should be projected back to the original dimension to construct the desired rotation matrix. All the constituent vectors in \(\bm{Q}\) should be mutually orthogonal, which equivalently means that the generated vectors should be in the null space of all previous vectors. In 2 dimensions, \(\bm{Q}\) can be written as,
\begin{equation} \label{eq:2Drot}
    \bm{Q} = \begin{bmatrix} \cos(\theta) & -\sin(\theta) \\
                        \sin(\theta) & ~~\cos(\theta) \end{bmatrix},
\end{equation}
and in 3 dimensions, \(\bm{Q}\) can be defined using any of the Euler sequence rotations \cite{schaub_analytical_2018}. MATLAB and Python codes are provided in a GitHub repository 
{here}\footnote{\url{https://github.com/saeidtafazzol/positive_definite_parameterization}}. 

\begin{algorithm}
\caption{\(N\)-dimensional Rotation Matrix 
Parameterization}
\label{alg: rotation}
\begin{algorithmic}[1]
\STATE \textbf{Input}: Angles $\theta_{ij}$ for $i = 1 ,\ldots , N, \; j = 1, \ldots, N-i$ 
\STATE \textbf{Output:} an Orthonormal Matrix $\bm{Q}$
\STATE \textbf{Initialize:} $\bm{Q} \leftarrow \{\}$ \COMMENT{The set of orthonormal vectors}
\FOR{$i = 1$ \textbf{to} $N$}
    \STATE $S = \text{null-space}(\bm{Q}^{\top})$
    \STATE $B = \text{Basis} (S)$ \COMMENT{Gram–Schmidt process}
    \STATE Construct $\bm{v}_i^p$ using $\theta_{i \, (1, \ldots, N-i)}$ and Eq.~\eqref{eq:unit_vec}
    \STATE $\bm{v}_i = B \bm{v}_i^p$ \COMMENT{Project back to $N$-dimensional space}
    \STATE $\bm{Q} \leftarrow [\bm{Q},  \bm{v}_i]$ \COMMENT{Append $\bm{v}_i$ as a new column to $\bm{Q}$}
\ENDFOR
\end{algorithmic}
\end{algorithm}

\section{Zermelo's Navigation Problem}\label{sec: zermelo}

\begin{figure}[h!]
    \centering
    \includegraphics[width=0.7\linewidth]{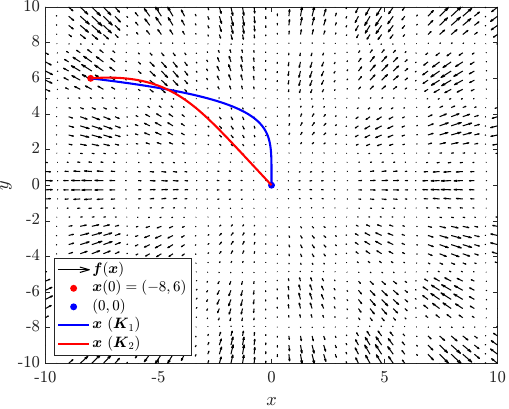}
    \caption{Zermelo's problem: state space with the trajectories from each solution with the vector field of nonlinearities in the dynamics.}
    \label{fig: zermelo phase space}
\end{figure}

\begin{figure}[h!]
    \centering
    \begin{subfigure}[b]{0.49\textwidth}
        \centering
        \includegraphics[width=\textwidth]{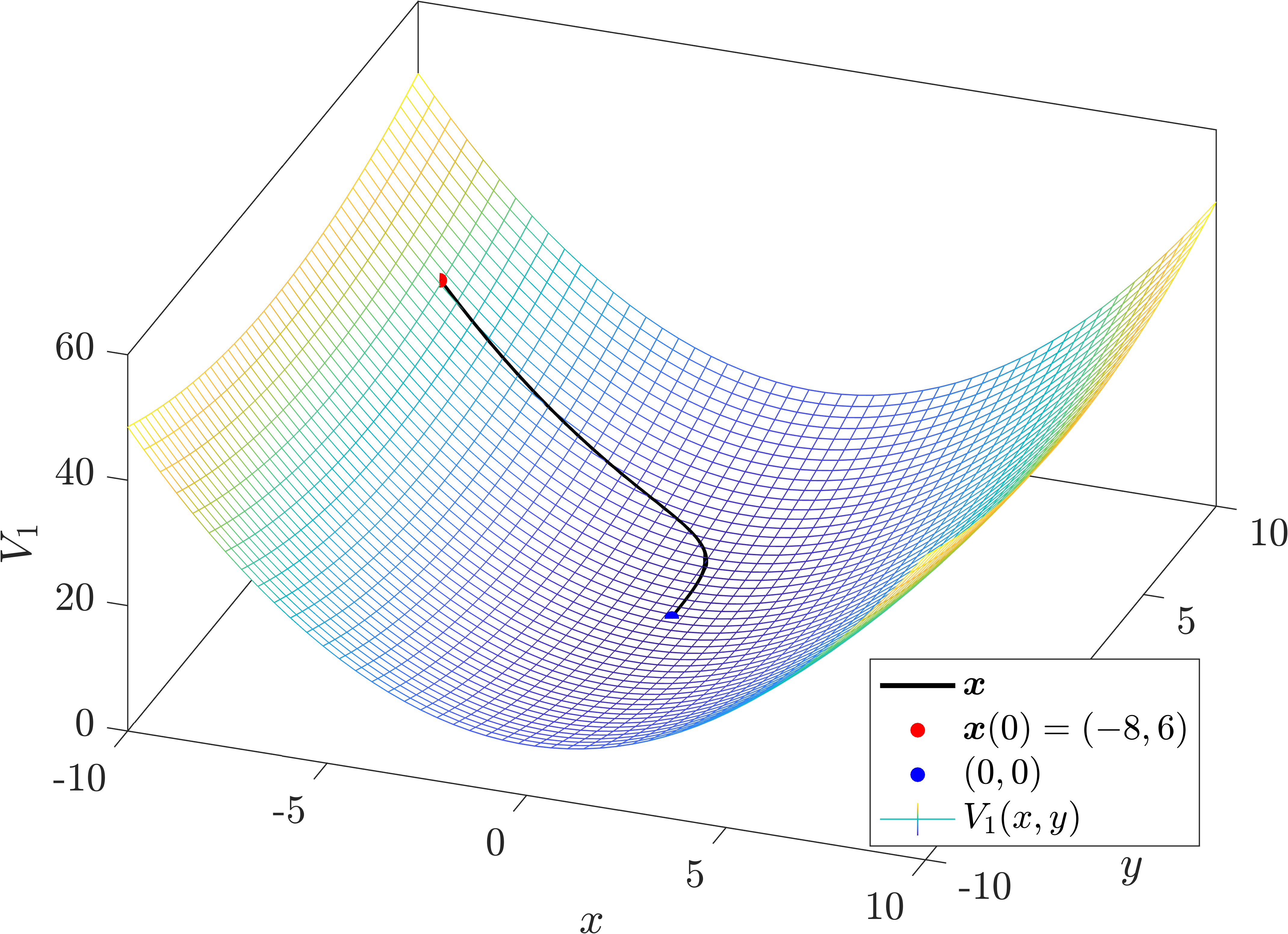}
        \caption{$\bm{K}_1$}
        \label{fig: zermelo clf 3d k1}
    \end{subfigure}
    \hfill
    \begin{subfigure}[b]{0.49\textwidth}
        \centering
        \includegraphics[width=\textwidth]{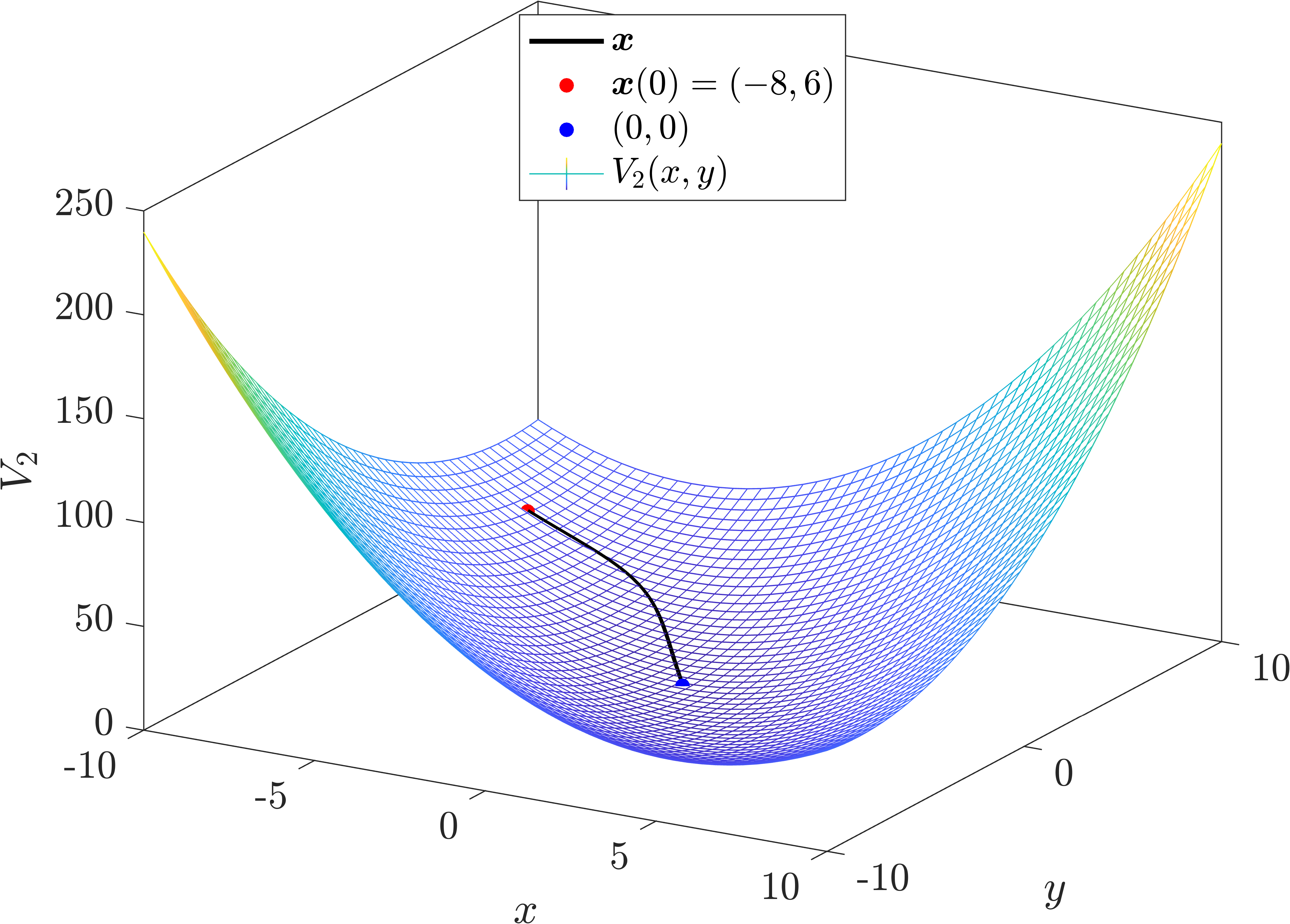}
        \caption{$\bm{K}_2$}
        \label{fig: zermelo clf 3d k2}
    \end{subfigure}
    \caption{Zermelo's problem: Lyapunov functions vs. states.}
    \label{fig: zermelo clf 3d}
\end{figure}

\begin{figure}[h!]
    \centering
    \begin{subfigure}[b]{0.49\textwidth}
        \centering
        \includegraphics[width=\textwidth]{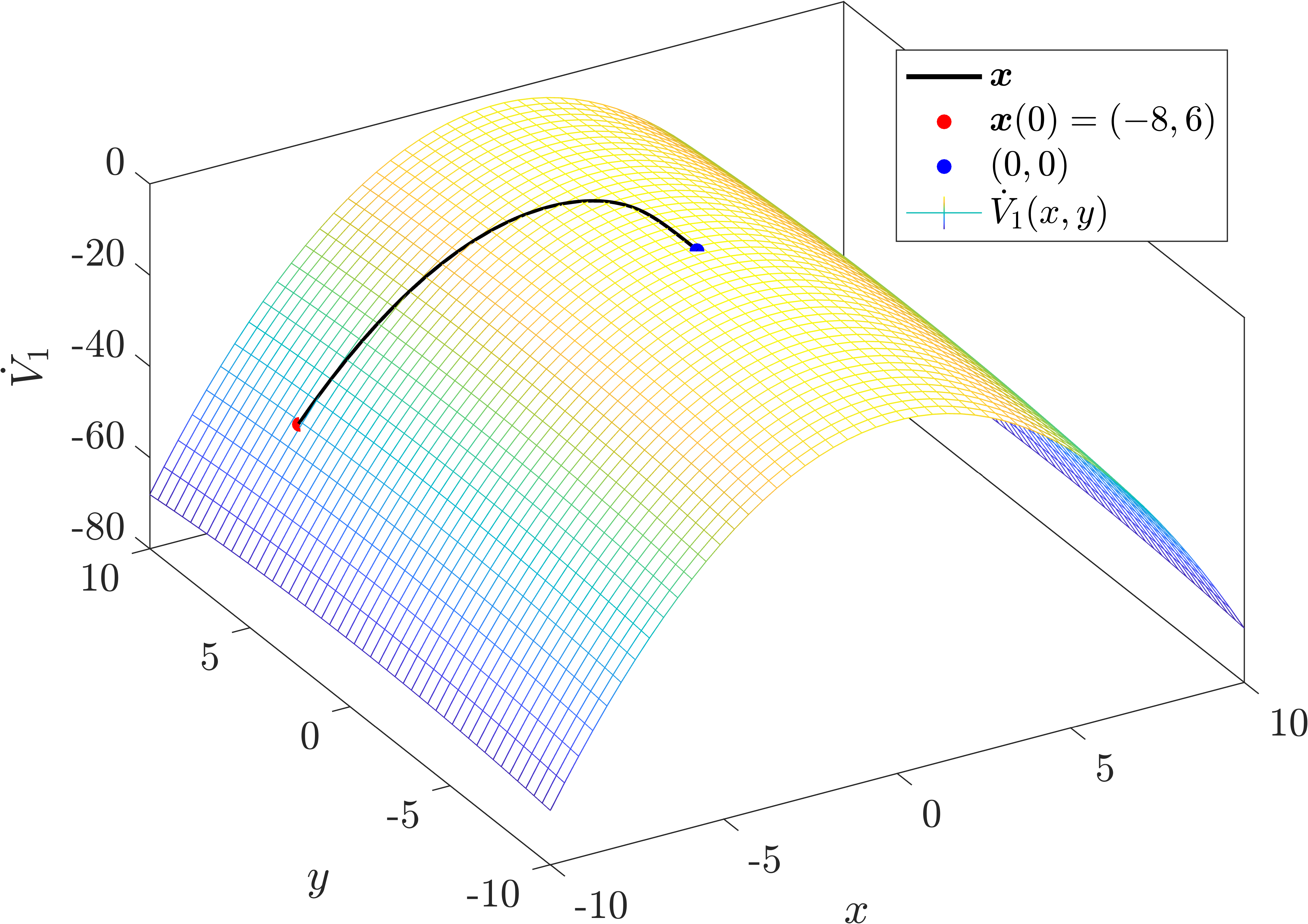}
        \caption{$\bm{K}_1$}
        \label{fig: zermelo clf time derivative 3d k1}
    \end{subfigure}
    \hfill
    \begin{subfigure}[b]{0.49\textwidth}
        \centering
        \includegraphics[width=\textwidth]{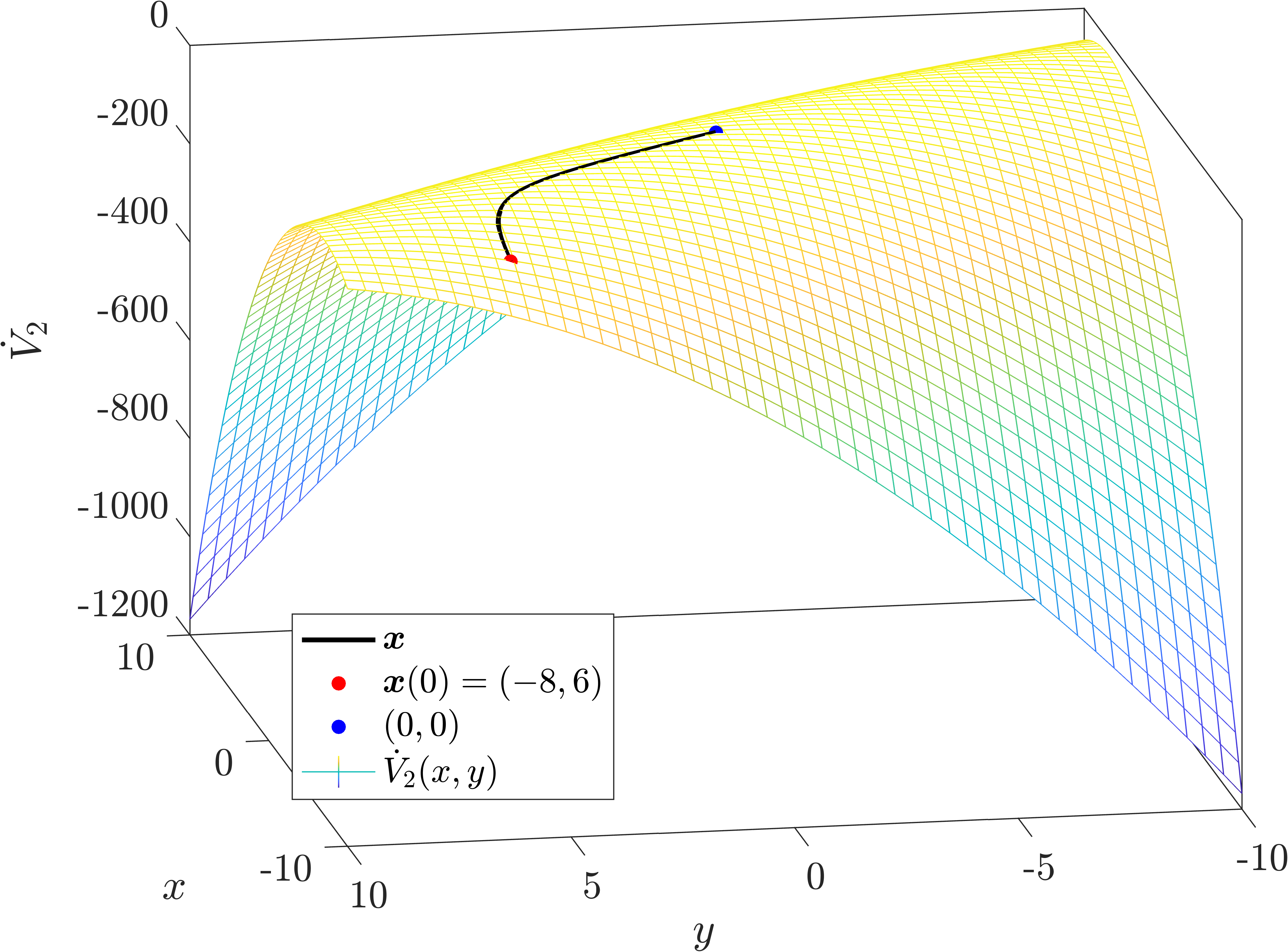}
        \caption{$\bm{K}_2$}
        \label{fig: zermelo clf time derivative 3d k2}
    \end{subfigure}
    \caption{Zermelo's problem: Lyapunov functions time-derivatives vs. states.}
    \label{fig: zermelo clf time derivative 3d}
\end{figure}

\begin{figure}[h!]
    \centering
    \begin{subfigure}[b]{0.49\textwidth}
        \centering
        \includegraphics[width=\textwidth]{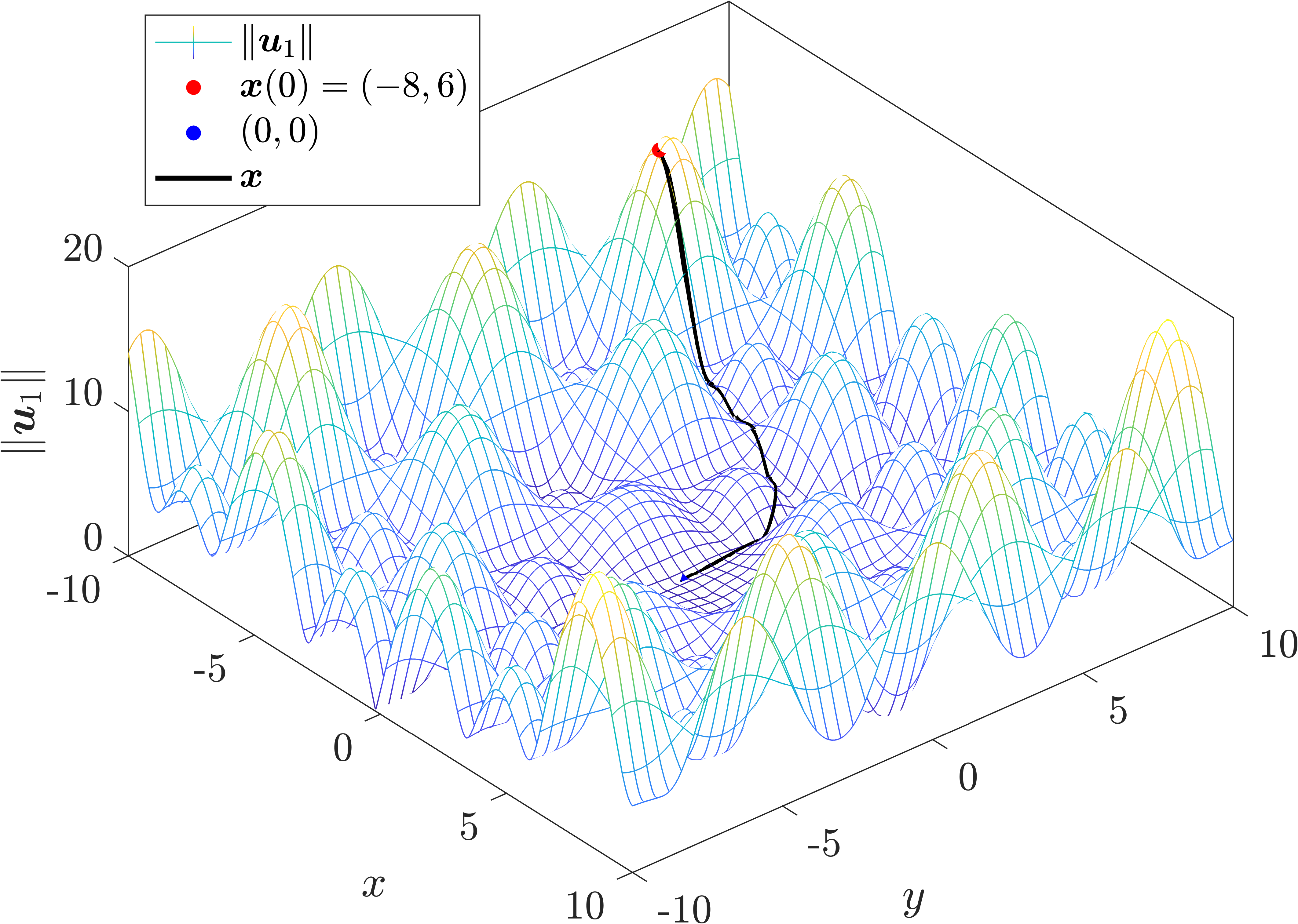}
        \caption{$\bm{K}_1$}
        \label{fig: zermelo control surface k1}
    \end{subfigure}
    \hfill
    \begin{subfigure}[b]{0.49\textwidth}
        \centering
        \includegraphics[width=\textwidth]{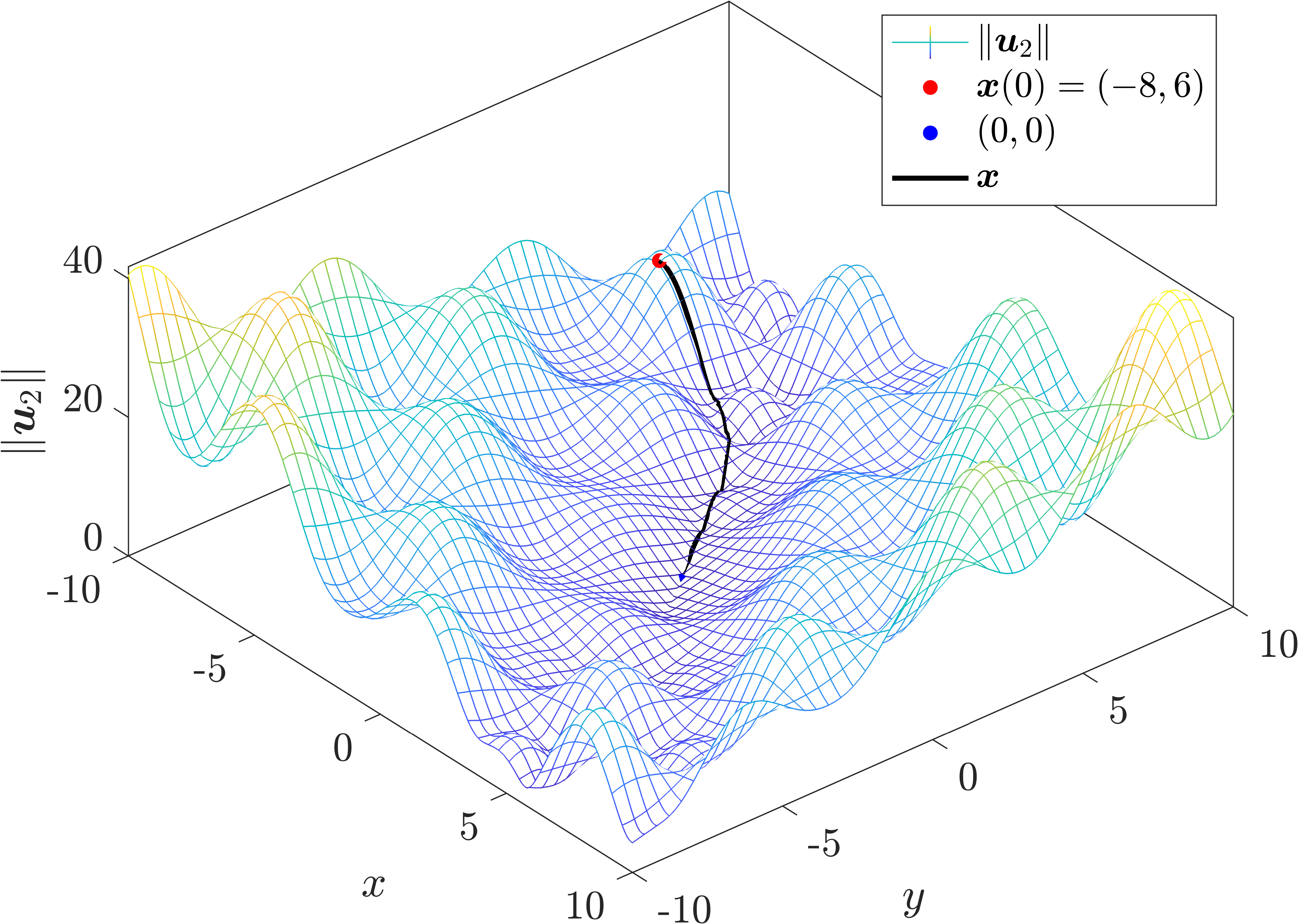}
        \caption{$\bm{K}_2$}
        \label{fig: zermelo control surface k2}
    \end{subfigure}
    \caption{Zermelo's problem: Euclidean norm of control function vs. states.}
    \label{fig: zermelo control surface}
\end{figure}

We solve a variation of Zermelo's navigation problem with a Lyapunov-based control law. The low-dimensionality of the problem allows for visualization of the candidate quadratic Lyapunov functions and the impact of using a full penalty matrix. Letting the states be \(\bm{x}^\top=[x,y]\) and the control be \(\bm{u}^\top=[u_x,u_y]\), the goal is to drive the system from an initial nonzero state, \(\bm{x}(0)\), to the origin, \([0,0]\) (i.e., the regularization problem). The dynamics are defined as,
\begin{align} \label{eq:ZermeloDynamics}
    \dot{\bm{x}} & = \bm{u} + \bm{f}(\bm{x}), & \bm{f}(\bm{x}) & =-[x,y]^{\top}\cos{y}\sin{x}.
\end{align}
The control, \(\bm{u}\), will be derived from a candidate quadratic Lyapunov function, which is written as,
\begin{equation} \label{eq: zermelo clf}
    V = \frac{1}{2}\bm{x}^\top\bm{K}\bm{x},
\end{equation}
where \(\bm{K}\in\mathbb{R}^{2\times2}\) is the positive-definite penalty matrix. The control law is derived such that the time derivative of Eq.~\eqref{eq: zermelo clf}, 
\begin{equation}
    \dot{V}=\frac{\partial V}{\partial \bm{x}}\dot{\bm{x}}=\bm{x}^\top\bm{K}(\bm{u} + \bm{f}(\bm{x})),
\end{equation}
is negative definite, which can be achieved by considering a control $\bm{u} = -\bm{f}(\bm{x}) - (\bm{x}^\top\bm{K})^\top$. Let \(\bm{K}_1=\text{diag}(k_1,k_2)\) denote a diagonal penalty matrix. We parameterize the full penalty matrix with eigendecomposition, \(\bm{K}_2\), as,
\begin{equation}
    \bm{K}_2 = \bm{Q}\bm{\Lambda}\bm{Q}^\top,
\end{equation}
where \(\bm{\Lambda}=\text{diag}(\lambda_1,\lambda_2)\) denotes the eigenvalue matrix and \(\bm{Q}\) is parametrized using Eq.~\eqref{eq:2Drot}. PSO is used to optimize the parameters (\(k_1 \in [0,10]\), \(k_2 \in [0,10]\) and \(\theta \in [0,2\pi]\)) to solve the optimization problem given in Eq.~\eqref{eq: min energy cost}. 
\begin{equation} \label{eq: min energy cost}
    \min_{\bm{u}}\quad J=\frac{1}{2}\int_{0}^\infty \bm{u}^\top\bm{u}~dt.
\end{equation}
We consider \(\bm{x}(0) = [-8,6]^\top\) and the dynamics in Eq.~\eqref{eq:ZermeloDynamics} are integrated over a time horizon of 100 seconds. An event-detection feature is used to terminate integration when \(\|\bm{x}\| \leq 10^{-3}\). The full and diagonal penalty matrices resulted in \(J=24.83\) and \(J=26.52\), respectively. The penalty matrices for each parameterization method were
\begin{align*}
    \bm{K}_1 = \begin{bmatrix} 0.8094 & 0 \\ 0 & 0.1611 \end{bmatrix},~ 
    \bm{K}_2 = \begin{bmatrix} 1.7421 & 0.9560 \\ 0.9560 & 1.1414 \end{bmatrix}.
\end{align*}
Figure \ref{fig: zermelo phase space} shows the state space with trajectories from both solutions, along with the nonlinearities in the dynamics shown as a vector field. Figures \ref{fig: zermelo clf 3d} and \ref{fig: zermelo clf time derivative 3d} show the CLFs and their time derivatives, respectively, plotted as a function of states with the trajectory plotted on top. Additionally, Figure \ref{fig: zermelo control surface} shows the surface of the control Euclidean norm as a function of states with the state trajectory overlaid. The difference in the CLFs, CLF time derivative, and control norm surfaces as a function of states evidently provide some advantage leading to a reduced control-effort requirement.

\section{Spacecraft Attitude Control Problem} \label{sec: spacecraft attitude control}

We proceed to compare the benefits of full penalty matrices in designing LQR and Lyapunov-based control laws for two arbitrarily chosen spacecraft attitude maneuvers: detumbling and rest-to-rest. The target state for both maneuvers is \((\psi(\infty),\theta(\infty),\phi(\infty))=(0,0,0)\) deg and 0 deg/s angular velocity (note that \(t=\infty\) is considered as the final time since the maneuvers are solved as infinite-time horizon problems). Both maneuvers start from an orientation of \((\psi(t_0),\theta(t_0),\phi(t_0))=(60,80,-60)\) deg and the detumbling maneuver has an initial angular velocity of \((p(t_0),q(t_0),r(t_0))=(0.1,-0.1,0.1)\) rad/s. The spacecraft has a moment of inertia matrix of \(\bm{I}=\text{diag}(10, 15, 20)\) kg\(\cdot\)m\(^2\).  These problem parameters were arbitrarily chosen and the coordinate representations (i.e., Euler angles and quaternions) are typically used in applications \cite{schaub_analytical_2018,junkins_optimal_1986}.

\subsection{LQR Attitude Control Design}

The spacecraft's orientation is parameterized using Euler angles, \(\psi\), \(\theta\), and \(\phi\) with angular velocity vector as \(\bm{\omega}^\top=\left[p,q,r\right]\). Let \(\bm{X}^\top = \begin{bmatrix} \psi ,\theta, \phi, \bm{\omega}^\top \end{bmatrix}\) denote the state vector. The attitude dynamics can be written as,
\begin{align}
    \dot{\bm{X}} & = \begin{bmatrix} \frac{\bm{\omega}}{c{\left(\theta\right)}} \begin{bmatrix} 
    0 & s{\left(\phi\right)} & c{\left(\phi\right)} \\ 
    0 & c{\left(\phi\right)}c{\left(\theta\right)} &  -s{\left(\phi\right)}c{\left(\theta\right)} \\ 
    c{\left(\theta\right)} & s{\left(\phi\right)}s{\left(\theta\right)} & c{\left(\phi\right)}s{\left(\theta\right)}\end{bmatrix}  \\ \bm{I}^{-1}\left(\bm{U} - \bm{\omega} \times \bm{I}\bm{\omega}\right)
    \end{bmatrix},
\end{align}
where $c(.) = \cos(.)$ and $s(.) = \sin(.)$. We linearize dynamics around the reference state and control vectors, \(\bm{X}_0=\bm{0}\) and \(\bm{U}_0=\bm{0}\), respectively. The state and control errors are denoted by \(\bm{x} = \bm{X}-\bm{X}_0\) and \(\bm{u}= \bm{U}-\bm{U}_0\), respectively. The quadratic cost functional whose minimization gives the optimal controller gain, i.e., \(\bm{k}_\text{LQR}\) appearing in \(\bm{U}\) (following the standard solution approach to LQR problems \cite{bryson_applied_1975}), is defined as,
\begin{align} \label{eq: lqr cost functional}
    \min_{\bm{k}_\text{LQR}}\quad J_\text{LQR} = \frac{1}{2}\int_{t_0}^\infty \left(\bm{x}^\top \bm{Q} \bm{x} + \bm{u}^\top \bm{R} \bm{u}\right) dt,
\end{align}
where \(\bm{Q}_{6\times6}\) and  \(\bm{R}_{3\times3}\) denote the penalty matrices. The control law used for propagating the nonlinear dynamics is \(\bm{U} = \bm{u} = - \bm{k}_\text{LQR} \bm{x}\) and the MATLAB function \verb|lqr| is used to find the optimal gain matrix, \(\bm{k}\). 
\begin{remark}\label{remark on N}
    There may be states that we do not want to be penalized. When \(\bm{Q}\) is diagonal, we can simply set the element corresponding to the state we do not want penalized to be 0. However, this is not as trivial when \(\bm{Q}\) is full. Rather than reducing the dimension of \(\bm{x}\) in Eq.~\eqref{eq: lqr cost functional} and thus possibly having to recode the LQR solution back into existing software architecture, a diagonal matrix \(\bm{N}\) with binary elements can be premultiplied with \(\bm{x}\) instead to choose which states are included in the cost functional. For an example, this is done with Q-law in Section \ref{sec: low thrust trajectory optimization}.
\end{remark}

PSO is used to tune the penalty matrices such that the maneuver is performed with minimal control effort (i.e., the same cost function given in Eq.~\eqref{eq: min energy cost}). 
Note that control takes the same form, \(\bm{u} = - \bm{k} \bm{x}\), but is applied on the state, $\bm{x}$, that is obtained from the nonlinear dynamics. Other applications may call for minimizing the deviation from a reference trajectory, in which case the cost function that PSO minimizes might be state-dependent. The diagonal penalty matrices are denoted as \(\bm{Q}_1\) and \(\bm{R}_1\) and the full penalty matrices are denoted as \(\bm{Q}_2\) and \(\bm{R}_2\). The upper and lower bounds for the diagonal elements of \(\bm{Q}_1\) and \(\bm{R}_1\) and for the eigenvalues of \(\bm{Q}_2\) and \(\bm{R}_2\) were 10 and \(10^{-8}\), respectively. 

PSO was invoked 10 times for each case to account for stochastic variations across solutions. Default options were used except for a swarm size of 500 and a maximum number of iterations of \(10^4\). These values were arbitrarily chosen such that solutions we deemed to be optimal for each respective form of penalty matrix could be obtained after a number of trials. While the theoretical time horizon considered is infinite, the maximum integration time is set to 100 seconds. The spacecraft attitude dynamics are propagated with MATLAB's \verb|ode113| with \verb|AbsTol| and \verb|RelTol| of \(1.0 \times 10^{-8}\). The built-in event-detection feature was used to determine when the states were ``close enough'' to the target state, \(\bm{X}_0\). The considered minimum-energy cost leads to solutions that use the entire time horizon (i.e., less aggressive control commands are prioritized).

\begin{figure}[h!]
    \centering
   \includegraphics[width=0.7\linewidth]{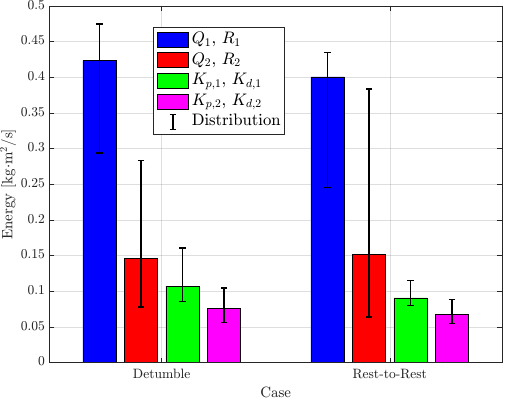}
    \caption{Comparison of the cost values for the LQR and Lyapunov-based controller with diagonal and full penalty matrices.}
    \label{fig: attitude control bar graph}
\end{figure}

\begin{figure}[h!]
    \centering
    \begin{subfigure}[b]{0.3\textwidth}
        \centering
        \includegraphics[width=\textwidth]{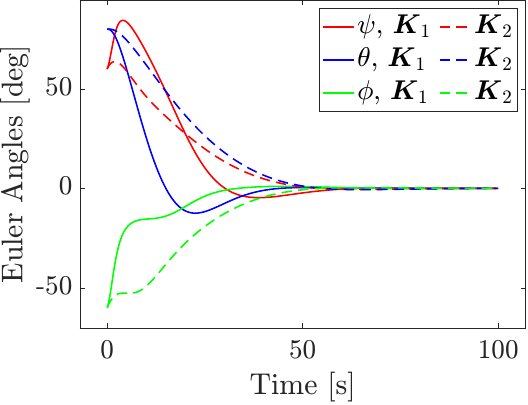}
        \caption{Euler angles vs. time.}
        \label{fig: lqr orientation manuever 1}
    \end{subfigure}
    \hfill
    \begin{subfigure}[b]{0.3\textwidth}
        \centering
        \includegraphics[width=\textwidth]{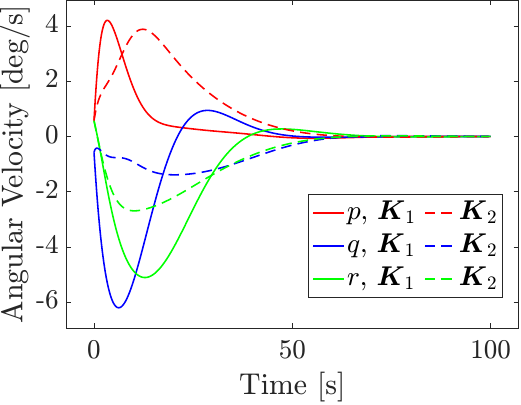}
        \caption{Angular velocity vs. time.}
        \label{fig: lqr angular velocity manuever 1}
    \end{subfigure}
    \hfill
    \begin{subfigure}[b]{0.3\textwidth}
        \centering
        \includegraphics[width=\textwidth]{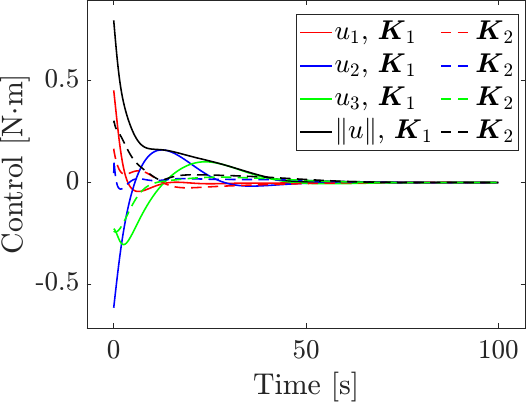}
        \caption{Control vs. time.}
        \label{fig: lqr control manuever 1}
    \end{subfigure}
\caption{Results from best detumbling solution for each type of penalty matrix with the LQR controller.}
\label{fig: lqr manuever 1 plots}
\end{figure}

Table \ref{tab: lqr results} summarizes the results for both maneuvers. Figure \ref{fig: attitude control bar graph} summarizes these results in a bar graph. A significant improvement in using the full penalty matrices was observed, with an average 65.4697\% and 62.0296\% decrease in cost value for the detumbling and rest-to-rest maneuvers, respectively. Figure \ref{fig: lqr manuever 1 plots} shows the most optimal solutions to the detumbling maneuver from both types of penalty matrices. The penalty matrices are given in \ref{app1}. 

\begin{table}[h!]
    \centering
    \caption{Detumbling and rest-to-rest maneuver results for LQR controllers with diagonal and full penalty matrices.}
    \label{tab: lqr results}
    \renewcommand{\arraystretch}{1.2} 
    \begin{tabularx}{\linewidth}{>{\centering\arraybackslash}p{1.5cm}>{\centering\arraybackslash}X>{\centering\arraybackslash}X>{\centering\arraybackslash}X>{\centering\arraybackslash}X}
        \hline 
         & \multicolumn{4}{c}{Minimum-Energy Cost [kg\(\cdot\)m\(^2\)/s]} \\ \hline
         & \multicolumn{2}{c}{Detumbling} & \multicolumn{2}{c}{Rest-to-rest} \\ \hline
        Run \# & \(\bm{Q}_1\), \(\bm{R}_1\) & \(\bm{Q}_2\), \(\bm{R}_2\) & \(\bm{Q}_1\), \(\bm{R}_1\)& \(\bm{Q}_2\), \(\bm{R}_2\)\\ \hline 
        1 & 0.47437 & 0.14191 & 0.43308 & 0.06817 \\ \hline 
        2 & 0.47436 & 0.09139 & 0.43307 & 0.06420 \\ \hline 
        3 & 0.30796 & 0.14311 & 0.29337 & 0.20775 \\ \hline 
        4 & 0.29380 & 0.08996 & 0.43307 & 0.09513 \\ \hline 
        5 & 0.47437 & 0.08397 & 0.43308 & 0.07716 \\ \hline 
        6 & 0.47437 & 0.25393 & 0.24552 & 0.11512 \\ \hline 
        7 & 0.47437 & 0.10233 & 0.43308 & 0.18937 \\ \hline 
        8 & 0.30785 & 0.28289 & 0.43424 & 0.07248 \\ \hline 
        9 & 0.47436 & 0.19370 & 0.43309 & 0.38302 \\ \hline 
        10 & 0.47435 & 0.07756 & 0.43307 & 0.24820 \\ \hline \hline
        Mean & 0.42302& 0.14607& 0.40047& 0.15206 \\ \hline 
        Best & 0.29380 & 0.07756 & 0.24552 & 0.06420 \\ \hline 
    \end{tabularx}
\end{table}

\subsection{Lyapunov-Based Attitude Control Design}

The Lyapunov-based controller uses quaternions, \(\bm{q}^\top=\left[q_0,q_1,q_2,q_3\right]\). This quaternion attitude representation is considered in this example to demonstrates the generalization of the method of the proposed and that in practice attitude parameterization is achieved with quaternions. The spacecraft's orientation rate is described in terms of the angular velocity vector \(\bm{\omega}^\top=\left[p,q,r\right]\). We propagate the spacecraft's attitude dynamics as
\begin{align}
    \bm{X} & = \begin{bmatrix} \bm{q} \\ \bm{\omega} \end{bmatrix}, & \dot{\bm{X}} & = \begin{bmatrix} \frac{1}{2}\bm{\Omega}\bm{q} \\ \bm{I}^{-1}\left(\bm{U} - \bm{\omega} \times \bm{I}\bm{\omega}\right)
    \end{bmatrix},
\end{align}
where \(\bm{\Omega}\) is a \(4\times4\) skew-symmetric form of $\bm{\omega}$ \cite{schaub_analytical_2018}. Using a nonlinear candidate quadratic Lyapunov function, the controller is derived in Ref. \cite{schaub_analytical_2018} and is defined as,
\begin{equation}
    \bm{U} = -K_p\bm{q}_\text{e} - K_d\bm{\omega},
\end{equation}
where $\bm{q}_\text{e} \in \mathbb{R}^3$ denotes the quaternion error vector and both \(K_p\) and \(K_d\) denoting \(3\times3\) control gain matrices (that are used in forming the Lyapunov function). Similar to the LQR problem, we denote the diagonal penalty matrices as \(\bm{K}_{p,1}\) and \(\bm{K}_{d,1}\) and the full penalty matrices as \(\bm{K}_{p,2}\) and \(\bm{K}_{d,2}\).

PSO is invoked 10 times for each case to account for stochastic variations across solutions and the penalty matrices are optimized to minimize the cost functional defined in Eq.~\eqref{eq: min energy cost}. The upper and lower bounds for the diagonal elements of \(\bm{K}_{p,1}\) and \(\bm{K}_{d,1}\) and for the eigenvalues of \(\bm{K}_{p,2}\) and \(\bm{K}_{d,2}\) were 10 and \(10^{-8}\). The maximum integration time is set to 100 seconds to create an artificial finite time horizon. The built-in event-detection feature of MATLAB's \verb|ode113| was used to determine when the states were ``close enough'' to the target states. The spacecraft attitude dynamics are propagated with MATLAB's \verb|ode113| with the same tolerances and PSO was run with the same options used in the LQR problem.

\begin{table}[h!]
    \centering
    \caption{Detumbling and rest-to-rest maneuver results for the Lyapunov-based controllers with diagonal and full penalty matrices.}
    \label{tab: lyapunov based attitude control results}
    \renewcommand{\arraystretch}{1.2} 
    \begin{tabularx}{\linewidth}{>{\centering\arraybackslash}p{1.5cm}>{\centering\arraybackslash}X>{\centering\arraybackslash}X>{\centering\arraybackslash}X>{\centering\arraybackslash}X}
        \hline 
         & \multicolumn{4}{c}{Minimum-Energy Cost [kg\(\cdot\)m\(^2\)/s]} \\ \hline
         & \multicolumn{2}{c}{Detumbling} & \multicolumn{2}{c}{Rest-to-rest} \\ \hline
        Run \# & \(\bm{K}_{p,1}\), \(\bm{K}_{d,1}\) & \(\bm{K}_{p,2}\), \(\bm{K}_{d,2}\) & \(\bm{K}_{p,1}\), \(\bm{K}_{d,1}\)& \(\bm{K}_{p,2}\), \(\bm{K}_{d,2}\)\\ \hline 
        1 & 0.09266 & 0.07841 & 0.08399 & 0.05507 \\ \hline 
        2 & 0.09665 & 0.10420 & 0.10394 & 0.05854 \\ \hline 
        3 & 0.09683 & 0.07002 & 0.08414 & 0.07439 \\ \hline 
        4 & 0.09362 & 0.09401 & 0.08073 & 0.07305 \\ \hline 
        5 & 0.10960 & 0.06194 & 0.10907 & 0.07109 \\ \hline 
        6 & 0.08541 & 0.05594 & 0.07993 & 0.08840 \\ \hline 
        7 & 0.08806 & 0.10104 & 0.08311 & 0.07015 \\ \hline 
        8 & 0.13964 & 0.05885 & 0.11482 & 0.05463 \\ \hline 
        9 & 0.16088 & 0.06611 & 0.08597 & 0.06427 \\ \hline 
        10 & 0.10820 & 0.06923 & 0.08034 & 0.06610 \\ \hline \hline
        Mean & 0.10716 & 0.07598& 0.09061& 0.06757 \\ \hline 
        Best & 0.08541 & 0.05594 &  0.07993 & 0.05463 \\
        \hline
        
    \end{tabularx}
\end{table}

\begin{figure}[h!]
    \centering
    \begin{subfigure}[b]{0.3\textwidth}
        \centering
        \includegraphics[width=\textwidth]{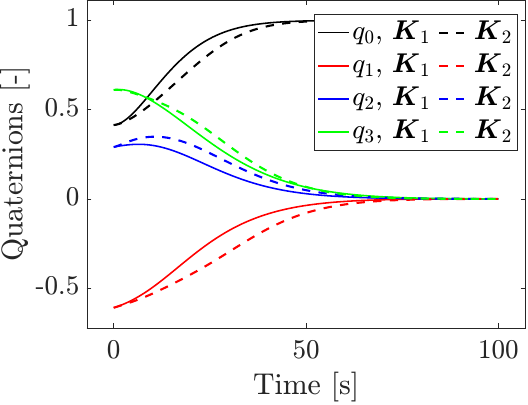}
        \caption{Quaternions vs. time.}
        \label{fig: lyapunov based orientation manuever 1}
    \end{subfigure}
    \hfill
    \begin{subfigure}[b]{0.3\textwidth}
        \centering
        \includegraphics[width=\textwidth]{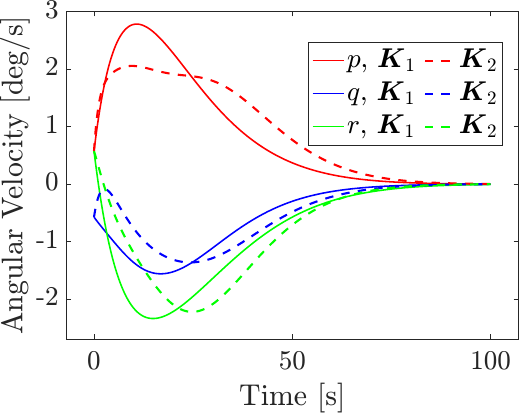}
        \caption{Angular velocity vs. time.}
        \label{fig: lyapunov based angular velocity manuever 1}
    \end{subfigure}
    \hfill
    \begin{subfigure}[b]{0.3\textwidth}
        \centering
        \includegraphics[width=\textwidth]{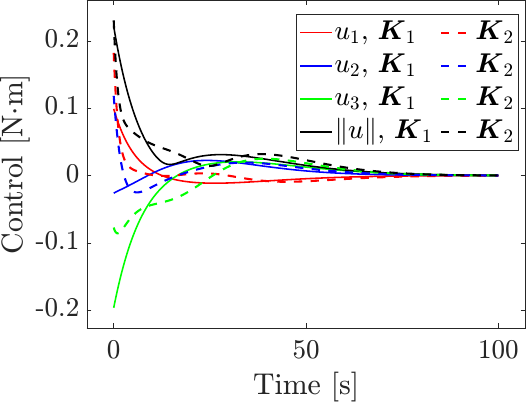}
        \caption{Control vs. time.}
        \label{fig: lyapunov based control manuever 1}
    \end{subfigure}
\caption{Results from best detumbling solution for each type of penalty matrix with the Lyapunov-based controller.}
\label{fig: lyapunov based manuever 1 plots}
\end{figure}


Table \ref{tab: lyapunov based attitude control results} shows the results for both maneuvers and with each type of penalty matrix. Figure \ref{fig: attitude control bar graph} shows a bar graph comparison of the values of costs. Again, improvement in using the full penalty matrices can be observed. However, the reduction in cost is not as significant as with the LQR controller. This difference is attributed to the significant improvement the nonlinear Lyapunov-based controller brings as opposed to the LQR controller which is based on dynamics linearized around the target state. Figure \ref{fig: lyapunov based manuever 1 plots} shows the solution for the most optimal solution to the detumbling maneuver from both types of penalty matrices. The penalty matrices are given in \ref{app1}.

\section{Low-Thrust Trajectory Optimization Problem} \label{sec: low thrust trajectory optimization}

To demonstrate the impact of using a full penalty matrix for Lyapunov-based low-thrust trajectory optimization, a generic LC law and the well-known Q-law \cite{petropoulos_low-thrust_2004} are used to solve a variety of low-thrust transfer trajectories. Both control laws are used with a diagonal penalty matrix, \(\bm{K}_1\), and with a full penalty matrix, \(\bm{K}_2\). The full penalty matrix, \(\bm{K}_2\), is generated using Algorithm \ref{alg: rotation}.

In the considered geocentric (and one Vesta-centered) low-thrust maneuvers, the goal is to transfer the spacecraft starting from a fully defined state to an orbit in minimum-time and, with Q-law only, in minimum-fuel. The boundary conditions and parameters for the five considered transfer scenarios are summarized in Table \ref{tab: orbit transfers}. These are a set of transfer cases originally defined in Ref.~\cite{petropoulos_low-thrust_2004} and have been solved as benchmark problems in a number of other studies, including Refs.~\cite{lee_design_2005} and \cite{varga_many-revolution_2016}. Note that the Case E* defines the boundary conditions that were used with Q-law. These set of boundary conditions, only differing in the inclination (as a result of a rotation of the inertial frame by \(30^\circ\) about the \(x\)-axis), was required to achieve convergence because we found that convergence was not possible with the implementation of Q-law as outlined in Ref.~\cite{petropoulos_low-thrust_2004} for the original Case E boundary conditions. 

\newgeometry{margin=2cm} 
\begin{landscape}
\begin{table}[]
    \footnotesize
    \centering
    \caption{Low-thrust orbit transfer cases.}
    \label{tab: orbit transfers}
    \renewcommand{\arraystretch}{1.4} 
    \begin{tabularx}{\linewidth}{>{\centering\arraybackslash}X>{\centering\arraybackslash}X>{\centering\arraybackslash}X>{\centering\arraybackslash}X>{\centering\arraybackslash}X>{\centering\arraybackslash}X>{\centering\arraybackslash}X>{\centering\arraybackslash}X>{\centering\arraybackslash}X>{\centering\arraybackslash}X>{\centering\arraybackslash}X>{\centering\arraybackslash}X>{\centering\arraybackslash}p{1.25cm}}
        \hline
        Cases & Orbit & \(a\) [deg] & \(e\) [-] & \(i\) [deg] & \(\Omega\) [deg] & \(\omega\) [deg] & Thrust [N] & Initial Mass [kg] & Specific Impulse [s] & Central Body & PSO Swarm Size & PSO Maximum Iterations \\ 
        \hline
        \multirow{2}{*}{A} & Initial & 7000 & 0.01 & 0.05 & 0 & 0 & \multirow{2}{*}{1} & \multirow{2}{*}{300} & \multirow{2}{*}{3100} & \multirow{2}{*}{Earth} & \multirow{2}{*}{50} & \multirow{2}{*}{50} \\ 
        \cline{2-7}
        & Target & 42000 & 0.01 & free & free & free & & & & & & \\ 
        \hline
        \multirow{2}{*}{B} & Initial & 24505.9 & 0.725 & 7.05 & 0 & 0 & \multirow{2}{*}{0.35} & \multirow{2}{*}{2000} & \multirow{2}{*}{2000} & \multirow{2}{*}{Earth} & \multirow{2}{*}{50} & \multirow{2}{*}{50} \\ 
        \cline{2-7}
        & Target & 42165 & 0.001 & 0.05 & free & free & & & & & & \\ 
        \hline
        \multirow{2}{*}{C} & Initial & 9222.7 & 0.2 & 0.573 & 0 & 0 & \multirow{2}{*}{9.3} & \multirow{2}{*}{300} & \multirow{2}{*}{3100} & \multirow{2}{*}{Earth} & \multirow{2}{*}{50} & \multirow{2}{*}{50} \\ 
        \cline{2-7}
        & Target & 30000 & 0.7 & free & free & free & & & & & & \\ 
        \hline
        \multirow{2}{*}{D} & Initial & 944.64 & 0.015 & 90.06 & -24.60 & 156.90 & \multirow{2}{*}{0.045} & \multirow{2}{*}{950} & \multirow{2}{*}{3045} & \multirow{2}{*}{Vesta} & \multirow{2}{*}{50} & \multirow{2}{*}{50} \\ 
        \cline{2-7}
        & Target & 401.72 & 0.012 & 90.01 & -40.73 & free & & & & & & \\ 
        \hline
         & \multirow{2}{*}{Initial} & 24505.9 & 0.725 & 0.06 & 0 & 0 & \multirow{4}{*}{2} & \multirow{4}{*}{2000} & \multirow{4}{*}{2000} & \multirow{4}{*}{Earth} & \multirow{4}{*}{200} & \multirow{4}{*}{300} \\ 
        E & & (24505.9) & (0.7) & (30.06) & (180) & (180) & & & & & & \\ 
        \cline{2-7}
        (E*) & \multirow{2}{*}{Target} & 26500 & 0.7 & 116 & 180 & 180 & & & & & & \\ 
        & & (26500) & (0.7) & (86) & (180) & (180) & & & & & & \\ 
        \hline
    \end{tabularx}
\end{table}
\end{landscape}
\restoregeometry

Earth and Vesta gravitational parameters are 398600.49 km\(^3\)/s\(^2\) and 17.8 km\(^3\)/s\(^2\), respectively. A canonical scaling is performed on all Cartesian states, such that \(\mu\) (the central body's gravitational parameter) is 1 DU\(^3\)/TU\(^2\). For the geocentric problems, the distance unit (DU) is 6378.1366 km and the time unit (TU) is 806.8110 seconds. For the Vesta-centered problem, DU is 289 km and TU is 1164.4927 seconds.

The state of the spacecraft, \(\bm{x} \in \mathbb{R}^7\), consists of Cartesian position \(\bm{r}\) and velocity \(\bm{v}\)  vectors, defined in an inertial frame centered at the central body. All acceleration vectors are expressed with respect to the inertial frame. An additional state, \(m\), is used to track the spacecraft's mass. The spacecraft equations of motion are written as,
\begin{align} \label{eq: eom}
    \bm{x} & = \begin{bmatrix} \bm{r} \\ \bm{v} \\ m \end{bmatrix}, & \dot{\bm{x}} & = \begin{bmatrix} \bm{v} \\ -\frac{\mu}{r^3} \bm{r} + \frac{T_\text{max}}{m}\hat{\bm{\alpha}}\delta \\ -\frac{T_\text{max}}{I_\text{sp}g_0}\delta \end{bmatrix},
\end{align}
where \(r=\|\bm{r}\|\) is the spacecraft's distance from the central body, \(T_\text{max}\) is the spacecraft's maximum thrust, \(I_\text{sp}\) is the spacecraft's specific impulse, \(g_0\) is the acceleration of gravity on Earth at sea level defined as 9.80665 m/s\(^2\), and \(\hat{\bm{\alpha}}\) is the thrust steering unit vector (i.e., \(\hat{\bm{\alpha}}^\top\hat{\bm{\alpha}}=1\)) assumed to freely orient in space. Lastly, \(\delta \in [0,1]\) is the thruster throttle magnitude control, which for the minimum-time transfers is 1 during the entire transfer, but is allowed to vary between 0 and 1 for the minimum-fuel transfers solved with Q-law. 

The diagonal penalty matrix, \(\bm{K}_1\), is trivially constructed and the full penalty matrix, \(\bm{K}_2\), is calculated according to Algorithm \ref{alg: rotation}. MATLAB's PSO algorithm is used to optimize the parameters for each of the control laws. Because PSO is a stochastic meta-heuristic optimization algorithm, it is invoked 5 times for each control law and both penalty matrices. The swarm sizes and the number of iterations are summarized in Table \ref{tab: orbit transfers}. These values are  chosen such that the resulting solutions appear to be close to the optimal one. The diagonal elements of \(\bm{K}_1\) and the eigenvalues of \(\bm{K}_2\) are bounded arbitrarily between 0 and 100. 

\subsection{Generic Lyapunov Control Law}

Generic LC laws are defined based on classical orbital elements \cite{vallado_fundamentals_2022}. These results were originally presented in Ref. \cite{nurre_expanding_2024} (except for cost function values that were reported in terms of mass rather than time of flight). We define an error vector, \(\bm{w}\), based on each transfer case in Table \ref{tab: orbit transfers} to be used in deriving the control law. This error vector is defined such that it becomes $\bm{0}$ when the final boundary conditions in Table \ref{tab: orbit transfers} are satisfied. We define the error vector for Case E, \(\bm{w}_\text{E}\), first. This error vector is defined in terms of the specific angular momentum vector, \(\bm{h}^\top = [h_x, h_y, h_z]\), and the eccentricity vector, \(\bm{e}\). Definitions of each are as follows \cite{vallado_fundamentals_2022}:
\begin{align} 
    \bm{h} & = \bm{r}\times\bm{v}, & \bm{e} & = \frac{\left(v^2-\frac{\mu}{r}\right)\bm{r} - \left(\bm{r}^\top\bm{v}\right)\bm{v}}{\mu}.
\end{align}

The error vector for Case E, \(\bm{w}_\text{E}\), is defined as, 
\begin{equation}
    \bm{w}_\text{E} = \begin{bmatrix} \bm{h} - \bm{h}_\text{T} \\ \bm{e} - \bm{e}_\text{T} \end{bmatrix}\in\mathbb{R}^6,
\end{equation}
where the subscript `$\text{T}$' corresponds to the target orbit values. 
\begin{remark}
    The considered choice of the error vector is entirely arbitrary for each maneuver. In fact, for Case E, a vector of only 5 dimensions is needed since only 5 orbital elements are being targeted at the final boundary condition due to transfer-type of the considered maneuvers.
\end{remark}

For Cases A through D, the boundary conditions are defined in terms of a combination of the orbital elements consisting of the specific angular momentum magnitude \(h\), eccentricity \(e\), inclination \(i\), and right ascension of the ascending node \(\Omega\). Definitions of the orbital elements are as follows \cite{vallado_fundamentals_2022}:
\begin{align} 
    h & = \left\|\bm{h}\right\|, & e & = \left\|\bm{e}\right\|, \\ i & = \cos^{-1}{\left(\frac{h_z}{h}\right)}, & \Omega & = \cos^{-1}{\left(\frac{\hat{\bm{x}}\cdot\bm{n}}{\|\hat{\bm{x}}\|\|\bm{n}\|}\right)},
\end{align}
where \(\bm{n}=\hat{\bm{z}}\times\bm{h} = [n_x,n_y,n_z]^{\top}\) is the line of nodes vector and a quadrant check is performed on \(\Omega\) such that if \(n_{y}<0\) then \(\Omega = 360^\circ - \Omega\). 

The error vectors for Cases A through D are defined as, 
\begin{align}
    \bm{w}_\text{A} & = \begin{bmatrix} h - h_\text{T} \\ e - e_\text{T} \end{bmatrix}\in\mathbb{R}^2, & \bm{w}_\text{B} & = \begin{bmatrix} h - h_\text{T} \\ e - e_\text{T} \\ i - i_\text{T} \end{bmatrix}\in\mathbb{R}^3, \\ 
    \bm{w}_\text{C} & = \begin{bmatrix} h - h_\text{T} \\ e - e_\text{T}\end{bmatrix}\in\mathbb{R}^2, & \bm{w}_\text{D} & = \begin{bmatrix} h - h_\text{T} \\ e - e_\text{T} \\ i - i_\text{T} \\ \Omega - \Omega_\text{T} \end{bmatrix}\in\mathbb{R}^4.
\end{align}
The CLF for each of the transfer problems is defined as,
\begin{equation} \label{eq: CLF}
    V_{j,i} = \frac{1}{2}\bm{w}_j^\top\bm{K}_i\bm{w}_j,
\end{equation}
for \(i=\{1,2\}\) and \(j=\{\text{A},\text{B},\text{C},\text{D},\text{E}\}\). The control law is derived to make the total time derivative of Eq.~\eqref{eq: CLF} (note that indices \(j\) and \(i\) are not written for brevity),
\begin{align}
    \dot{V} & =\frac{\partial V}{\partial \bm{r}}\bm{v}+\frac{\partial V}{\partial \bm{v}}\left(-\frac{\mu}{r^3} \bm{r} + \frac{T_\text{max}}{m}\hat{\bm{\alpha}}\right),
\end{align}
becomes negative, leading to the control law, in Eq.~\eqref{eq: control law}, as,
\begin{equation} \label{eq: control law}
    \hat{\bm{\alpha}}^*=-\left(\frac{\partial V}{\partial \bm{v}}\right)^\top/\left\|\frac{\partial V}{\partial \bm{v}}\right\|.
\end{equation}
We used automatic differentiation capabilities of CasADi \cite{andersson_casadi_2019} to derive $\hat{\bm{\alpha}}^*$. MATLAB's variable-step variable-order explicit differential equation solver \verb|ode113| is used to integrate Eq.~\eqref{eq: eom} with the control law in Eq.~\eqref{eq: control law} substituted in for \(\hat{\bm{\alpha}}\) with absolute and relative integration tolerances of \(10^{-10}\). The event-detection capability of \verb|ode113| is used to determine when the solution is close enough to the target orbit, since LC laws will not make the system converge in finite time.

\newgeometry{margin=2cm} 
\begin{landscape}
\begin{table}[]
    \footnotesize
    \centering
    \caption{Low-thrust results for generic Lyapunov-based control laws.}
    \label{tab: generic control law low thrust results}
    \renewcommand{\arraystretch}{1.6} 
    \begin{tabularx}{\linewidth}{>{\centering\arraybackslash}m{1cm}>{\centering\arraybackslash}m{1.5cm}>{\centering\arraybackslash}X>{\centering\arraybackslash}X>{\centering\arraybackslash}X>{\centering\arraybackslash}X>{\centering\arraybackslash}X>{\centering\arraybackslash}X>{\centering\arraybackslash}X>{\centering\arraybackslash}m{3cm}}
        \hline
        \multirow{2}{*}{Cases} & \multirow{2}{=}{Penalty Matrix} & \multicolumn{6}{c}{Objective Value, \(J(K_{1,2})=t_f\) [days]} & \multirow{2}{=}{\(\text{mean}(J(K_2))-\text{mean}(J(K_1))\) [days]} \\ 
        \cline{3-8}
         & & Run 1 & Run 2 & Run 3 & Run 4 & Run 5 & Mean & \\ 
        \hline
        \multirow{2}{*}{A} & \(K_1\) & 14.5705 & 14.5702 & 14.57 & 14.5703 & 14.5701 & 14.5702 & \multirow{2}{*}{-0.0952252} \\ 
        \cline{2-8}
         & \(K_2\) & 14.4751 & 14.4748 & 14.4749 & 14.4754 & 14.4748 & 14.475 & \\ 
        \hline
        \multirow{2}{*}{B} & \(K_1\) & 142.229 & 142.229 & 142.229 & 142.229 & 142.229 & 142.229 & \multirow{2}{*}{-3.18558} \\ 
        \cline{2-8}
         & \(K_2\) & 139.03 & 139.067 & 139.02 & 139.074 & 139.023 & 139.043 & \\ 
        \hline
        \multirow{2}{*}{C} & \(K_1\) & 1.5102 & 1.5102 & 1.5102 & 1.5102 & 1.5102 & 1.5102 & \multirow{2}{*}{-0.018357} \\ 
        \cline{2-8}
         & \(K_2\) & 1.49184 & 1.49184 & 1.49184 & 1.49184 & 1.49184 & 1.49184 & \\ 
        \hline
        \multirow{2}{*}{D} & \(K_1\) & 24.9903 & 24.9903 & 24.9903 & 24.9903 & 24.9903 & 24.9903 & \multirow{2}{*}{-0.153639} \\ 
        \cline{2-8}
         & \(K_2\) & 24.6992 & 24.7059 & 24.9456 & 24.9456 & 24.8868 & 24.8366 & \\ 
        \hline
        \multirow{2}{*}{E} & \(K_1\) & 93.581 & 80.9959 & 83.4887 & 89.0044 & 106.926 & 90.7993 & \multirow{2}{*}{-2.2493} \\ 
        \cline{2-8}
         & \(K_2\) & 109.751 & 89.6992 & 85.8899 & 79.7213 & 77.6889 & 88.55 & \\ 
        \hline
    \end{tabularx}
\end{table}
\end{landscape}
\restoregeometry

Table \ref{tab: generic control law low thrust results} summarizes the results for both \(\bm{K}_1\) and \(\bm{K}_2\). It can be observed that considering \(\bm{K}_2\) in the CLF improves the time of flight in each transfer case, with the improvement being more evident for maneuvers occurring over a longer time horizon like Case B or that are more complicated like Case E. The global minimum is achievable for the shorter duration transfers, especially Cases A and C. This is substantiated by PSO being able to find roughly the same solution for all 5 runs of each penalty matrix. For these cases, the new parameterization of the penalty matrix improved the solution optimality given the improvement in the global minimum. 

\begin{figure}[]
    \centering
    \includegraphics[width=0.7\linewidth]{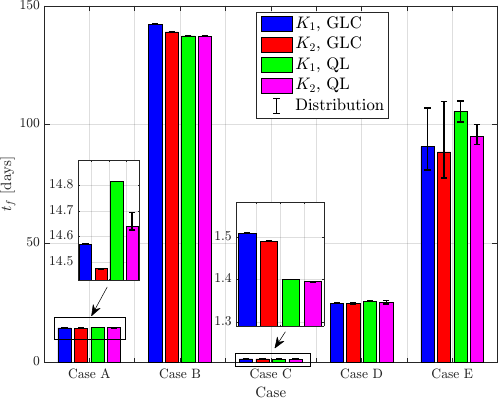}
    \caption{Average time-of-flight \(t_f\) for each case and penalty matrix for the low-thrust results for the generic Lyapunov-based control laws and Q-law, denoted ``GLC'' and ``QL'' in the legend, respectively.}
    \label{fig: min time low thrust bar graph}
\end{figure}

The results indicate that multiple local minima are found for the other transfer cases. In particular, Case E demonstrates the existence of many local minima, with the times of flight ranging broadly and with the worst solution achieved being found by the full penalty matrix, whereas in the other cases the worst solution was always found by the diagonal penalty matrix. On the other hand, the best solution is found by the full matrix representation in all transfer cases and the average times of flights for the full penalty matrix, \(\bm{K}_2\), are lower than those from the diagonal penalty matrix, \(\bm{K}_1\), in every transfer case. Figure \ref{fig: min time low thrust bar graph} summarizes the results in Table \ref{tab: generic control law low thrust results} as a bar graph. The average times of flight represent the bars and the error bars represent the distribution of time of flight for the 5 runs of PSO in each case. This error bar further illustrates the few or many local minima existing in each case, with Case E having the widest distribution.

\begin{figure}[]
    \centering
    \includegraphics[width=0.7\linewidth]{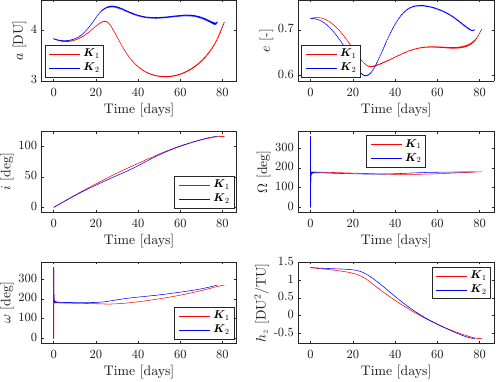}
    \caption{Case E orbital elements: minimum-time generic control law with \(\bm{K}_1\) and \(\bm{K}_2\).}
    \label{fig: glc oe}
\end{figure}

\begin{figure}[]
    \centering
    \begin{subfigure}[b]{0.49\textwidth}
        \centering
        \includegraphics[width=\textwidth]{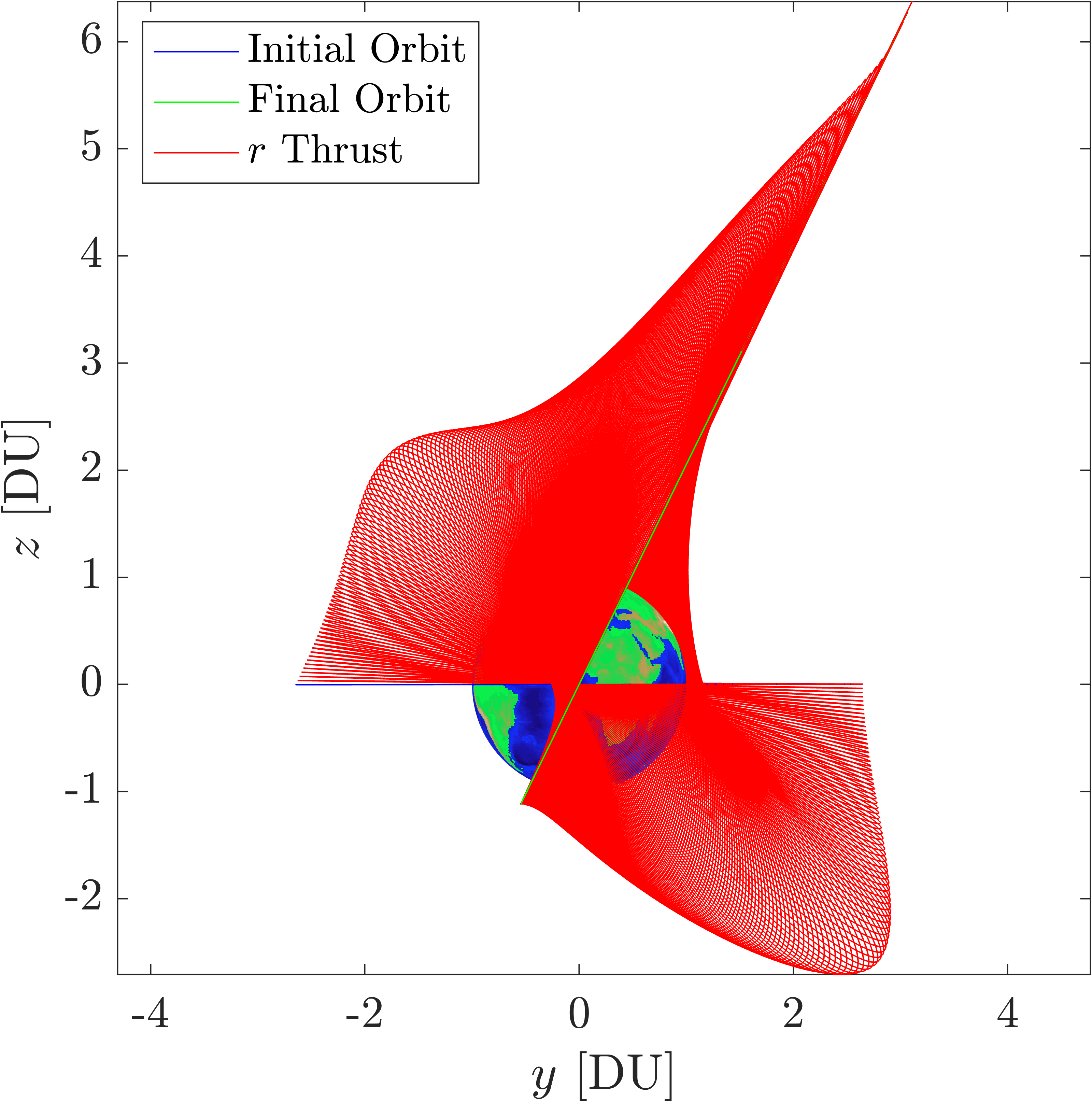}
        \caption{Case E trajectory: minimum-time generic control law with \(\bm{K}_1\).}
        \label{fig: glc k1 trajectory}
    \end{subfigure}
    \hfill
    \begin{subfigure}[b]{0.49\textwidth}
        \centering
        \includegraphics[width=\linewidth]{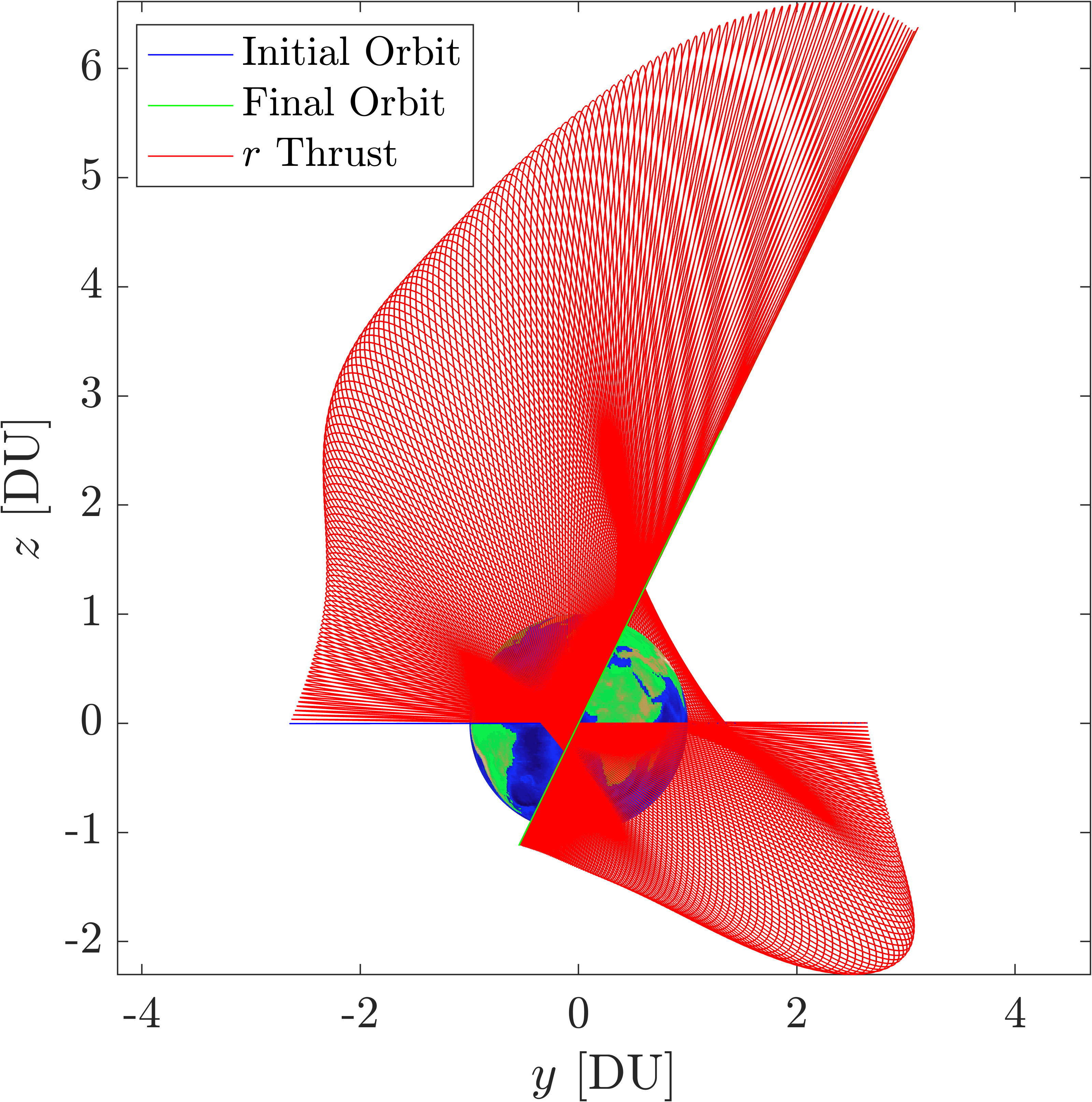}
        \caption{Case E trajectory: minimum-time generic control law with \(\bm{K}_2\).}
        \label{fig: glc k2 trajectory}
    \end{subfigure}
    \caption{Case E trajectory: minimum-fuel Q-law with \(\bm{K}_2\).}
    \label{fig: glc case e mf trajectory}
\end{figure}

Figure \ref{fig: glc oe} also shows the time histories of the orbital elements for the best Case E solution for \(\bm{K}_1\) and \(\bm{K}_2\). Note that the oscillations followed by a settling to a particular value in \(\Omega\) and \(\omega\) are due to the orbit being nearly equatorial at the beginning of the transfer, leading to \(\Omega\) and \(\omega\) not being properly defined. Figure \ref{fig: glc case e mf trajectory} shows the trajectory for these solutions. The full weighting matrix, \(\bm{K}_2\), drives the orbital elements to their target values differently (and likely more efficiently based on the improved times of flight) than the diagonal penalty matrix \(\bm{K}_1\). Inspecting the \(\bm{K}_1\) solution in red Figure \ref{fig: glc oe}, the semi-major axis has first increased to beyond the semi-major axis of the target orbit, but then, it is reduced noticeably below the target value. On the other hand, the eccentricity is first decreased and then gradually increased. Inspecting the \(\bm{K}_2\) solution in blue Figure \ref{fig: glc oe}, the semi major axis is instead increased and maintained more or less until the end of the maneuver. The eccentricity, however, has changed in a sinusoidal manner. The \(\bm{K}_2\) parameterization taking into account the cross-coupling errors of the orbital elements is likely what results in a solution with a better time of flight than the diagonal matrix \(\bm{K}_1\). 

\subsection{Q-law Low-Thrust Trajectory Design Method}
Q-law is used to solve the same cases in Table \ref{tab: orbit transfers}. The Q-law from Ref. \cite{petropoulos_low-thrust_2004} is used, but with some modifications to make it amenable to a full penalty matrix implementation. The proximity quotient with both \(\bm{K}_1\) and \(\bm{K}_2\), \(Q_{1,2}\), is defined as,
\begin{equation} \label{eq: Q}
    Q_{1,2}=\left(1+W_PP\right)\left(\bm{S}\bm{N}\bm{D}\right)^\top\bm{K}_{1,2}\left(\bm{N}\bm{D}\right),
\end{equation}
where \(\bm{D}\) denotes the error vector defined as,
\begin{align} \label{eq:ervecQlaw}
    \bm{D}&= \left(\frac{\text{\oe}-\text{\oe}_\text{T}}{\dot{\text{\oe}}_\text{xx}}\right), & \text{for} & ~\text{\oe}=a,e,i,\omega,\Omega.
\end{align}
In Eq.~\eqref{eq:ervecQlaw}, \(\dot{\text{\oe}}_\text{xx}\) denotes the maximum rate of change of the respective orbital element over thrust direction angles and over true anomaly on the osculating orbit. Analytical relations for these rates are given in Ref.~\cite{petropoulos_low-thrust_2004}. Note that the difference between instantaneous and target \(\omega\) and \(\Omega\) values is taken inside inverse cosine and cosine functions in Ref. \cite{petropoulos_low-thrust_2004}. We found that a singularity in the derivative of Eq.~\eqref{eq: Q} with respect to state (which is required to derive Q-law \cite{shannon_q-law_2020}) was encountered when the instantaneous element, \(\text{\oe}\), equals the target element, \(\text{\oe}_\text{T}\). For instance, for \(\omega\), we have
\begin{equation}
    \frac{d}{d\omega}\left(\cos^{-1}{\left(\cos{\left(\omega-\omega_\text{T}\right)}\right)}\right)=\frac{\sin{\left(\omega-\omega_\text{T}\right)}}{\sqrt{1-\cos^2{\left(\omega-\omega_\text{T}\right)}}}.
\end{equation}
We found for some transfer cases that this singularity can be encountered frequently before the transfer is complete.

The scaling matrix function, \(\bm{S}\) in Eq.~\eqref{eq: Q}, is defined as,
\begin{equation}
    \bm{S} = \text{diag}\left(S_a,S_e,S_i,S_\omega,S_\Omega\right),
\end{equation}
where the elements of \(\bm{S}\) are defined according to Ref. \cite{petropoulos_low-thrust_2004} as,
\begin{equation}
    S_\text{\oe} = \begin{cases}{\left[1+\left(\frac{a-a_T}{m~ a_T}\right)^n\right]^{\frac{1}{r}}}, & \text { for } \text{\oe}=a, \\ 1, & \text { for } \text{\oe}=e, i, \omega, \Omega,\end{cases}
\end{equation}
where \(m\), \(n\), and \(r\) are scalars with nominal values set to 3, 4, and 2, respectively. In Ref. \cite{petropoulos_low-thrust_2004}, if certain elements weren't desired to be targeted, then their penalty factor, \(W_\text{\oe}\), could simply be set equal to 0. When \(\bm{K}_1\) is used, the same step can be performed by setting the desired elements of \(\bm{K}_1\) to 0 to avoid targeting those elements. However, when \(\bm{K}_2\) is used, the matrix \(\bm{N}\) must be introduced and is defined as,
\begin{equation}
    \bm{N} = \text{diag}\left(N_a,N_e,N_i,N_\omega,N_\Omega\right),~\text{with}
\end{equation}
\begin{equation}
    N_\text{\oe} = \begin{cases} 1 & \text{if \oe~is targeted}\\ 0 & \text{if \oe~is not targeted} \end{cases} ~ \text{for}~\text{\oe}=a,e,i,\omega,\Omega.
\end{equation}
This follows the procedure described in Remark \ref{remark on N}, which can be useful when control laws must be hard-coded into vehicle hardware, or when it is cumbersome to rederive the control law depending on the problem being considered. Alternatively, Eq.~\eqref{eq: Q} could be rederived for a different set of \(\text{\oe}\) depending on what elements are desired to be targeted. However, the present approach allows for a more user-friendly implementation for solving a wider variety of problems.

The penalty function, \(P\), serves to enforce a minimum-periapsis-radius constraint and is defined as,
\begin{equation}
    P=\exp{\left(k\left(1-\frac{r_p}{r_{p,\text{min}}}\right)\right)},
\end{equation}
where the value of \(k\) is set to 4. The periapsis radius can be found with \(r_p=a(1-e^2)\) and \(r_{p,\text{min}}\) is the prescribed minimum-periapsis radius and is set arbitrarily to 10 km in this work for all problems. In Eq.~\eqref{eq: Q}, the penalty factor, \(W_P\), is either 1 or 0 to activate 0r deactivate the penalty function \(P\). In this work, it is set to 1.

The thrust steering control law, \(\hat{\bm{\alpha}}\), is parameterized by two angles, which are derived similarly to Ref. \cite{shannon_q-law_2020}. For minimum-time solutions, the thrust magnitude control \(\delta\) is always equal to 1. Ref.~\cite{petropoulos_refinements_2005} describes a coasting mechanism, which determines the thrust magnitude when minimum-fuel transfers are desired which is defined as,
\begin{equation}
    \delta=\begin{cases} 0 & \eta_r \leq \eta_\text{cut}, \\ 1 & \eta_r > \eta_\text{cut}, \end{cases}
\end{equation}
where \(\eta_\text{cut}\in[0,1]\) is a user-defined cut-off value, which determines the trade-off between time-of-flight and propellant savings and \(\eta_r\) is the so-called relative effectivity, which is a measure of how effective it is to thrust on an orbit and is defined as,
\begin{equation}
    \eta_r = \frac{\dot{Q}_\text{n}-\dot{Q}_\text{nx}}{\dot{Q}_\text{nn}-\dot{Q}_\text{nx}},
\end{equation}
where \(\dot{Q}_\text{n}\) is the time-derivative of \(Q\) minimized over the steering control, \(\dot{Q}_\text{nn}\) is \(\dot{Q}_\text{n}\) minimized over the orbit's true anomaly, and \(\dot{Q}_\text{nx}\) is \(\dot{Q}_\text{n}\) maximized over the orbit's true anomaly. Please refer to Ref.~\cite{petropoulos_refinements_2005} for further details. The spacecraft's dynamics are propagated in classical orbital elements (COEs) using MATLAB's \verb|ode45| with both absolute and relative tolerances set to \(10^{-10}\). This choice of numerical integrator was different from the one used in the previous sections due to existing code architecture of the authors. Checks are performed during integration to ensure that the singularities at circular and planar orbits are not encountered. Due to the chattering that can result from the Q-law \cite{noble_ariel_hatten_critical_2012,hecht_q-law_2024}, the control input is assumed constant over 1-minute intervals. 

\begin{figure}[h!]
    \centering
    \includegraphics[width=0.7\linewidth]{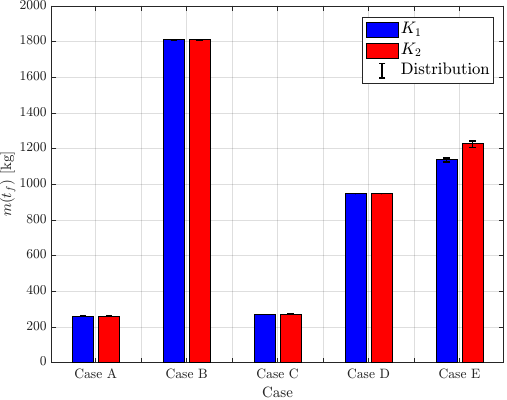}
    \caption{Average minimum-fuel objective values for the Q-law solutions.}
    \label{fig: qlaw min fuel bar graph}
\end{figure}

The results are summarized in Table \ref{tab: q law low thrust results}. The last column shows the difference in the mean objective value across the 5 runs for each penalty matrix. Because time-of-flight is being minimized and final fuel mass is being maximized, the full weighting matrix shows improvement in every case, with the difference being negative for minimum-time problems and positive for minimum-fuel problems. Most cases show little improvement in objective values, indicating that Q-law might be less dependent on the penalty matrix and more dependent on the maximum orbital element rates-of-change, \(\dot{\text{\oe}}_\text{xx}\). However, Case E shows a significantly larger improvement with a 10.5 day improvement in time-of-flight and a nearly 90 kg improvement in final fuel mass. Figure \ref{fig: min time low thrust bar graph} summarizes the minimum-time results (along with the minimum-time results from the generic Lyapunov control laws) where the average time-of-flight represents the bars and the error bars represent the distribution in solutions. The lower distribution visible in Case E further shows how Q-law might be less dependent on the value of the penalty matrix than the generic Lyapunov control law. Figure \ref{fig: qlaw min fuel bar graph} shows the results from the minimum-fuel solutions summarized as a bar graph.

\newgeometry{margin=2cm} 
\begin{landscape}
\begin{table}[]
    \footnotesize
    \centering
    \caption{Q-law low-thrust transfer results.}
    \label{tab: q law low thrust results}
    \renewcommand{\arraystretch}{1.6} 
    \begin{tabularx}{\linewidth}{>{\centering\arraybackslash}X>{\centering\arraybackslash}X>{\centering\arraybackslash}X>{\centering\arraybackslash}X>{\centering\arraybackslash}X>{\centering\arraybackslash}X>{\centering\arraybackslash}X>{\centering\arraybackslash}X>{\centering\arraybackslash}X>{\centering\arraybackslash}m{2.5cm}}
        \hline
        \multirow{2}{*}{Case} & \multirow{2}{=}{Objective, \(J(K_{1,2})\)} & \multirow{2}{=}{Penalty Matrix} & \multicolumn{6}{c}{Objective Value, \(J(K_{1,2})=t_f~\text{or}~m_f\) [days or kg]}  & \multirow{2}{=}{\(\text{mean}(J(K_2))-\text{mean}(J(K_1))\) [days or kg]} \\ 
        \cline{4-9}
         & & & Run 1 & Run 2 & Run 3 & Run 4 & Run 5 & Mean & \\ 
        \hline
        \multirow{4}{*}{A} & \multirow{2}{*}{\(t_f\) [days]} & \(K_1\) & 14.815 & 14.815 & 14.815 & 14.815 & 14.815 & 14.815 & \multirow{2}{*}{-0.175972} \\ \cline{3-9}
         & & \(K_2\) & 14.626 & 14.626 & 14.625 & 14.694 & 14.626 & 14.639 & \\ 
        \cline{2-10}
         & \multirow{2}{*}{\(m_f\) [kg]} & \(K_1\) & 259.807 & 259.807 & 259.809 & 259.809 & 259.807 & 259.808 & \multirow{2}{*}{0.496174} \\ 
        \cline{3-9}
         & & \(K_2\) & 260.304 & 260.346 & 260.306 & 260.239 & 260.324 & 260.304 & \\ 
        \hline
        \multirow{4}{*}{B} & \multirow{2}{*}{\(t_f\) [days]} & \(K_1\) & 137.31 & 137.309 & 137.31 & 137.31 & 137.31 & 137.31 & \multirow{2}{*}{-0.0606944} \\ \cline{3-9}
         & & \(K_2\) & 137.182 & 137.309 & 137.189 & 137.256 & 137.309 & 137.249 &  \\ 
        \cline{2-10}
         & \multirow{2}{*}{\(m_f\) [kg]} & \(K_1\) & 1809.67 & 1809.67 & 1809.66 & 1809.67 & 1809.71 & 1809.67 & \multirow{2}{*}{0.0409008} \\ 
        \cline{3-9}
         & & \(K_2\) & 1809.71 & 1809.67 & 1810 & 1809.73 & 1809.46 & 1809.72 &  \\ 
        \hline
        \multirow{4}{*}{C} & \multirow{2}{*}{\(t_f\) [days]} & \(K_1\) & 1.4007 & 1.4007 & 1.4007 & 1.4007 & 1.4007 & 1.4007 & \multirow{2}{*}{-0.0052778} \\ \cline{3-9}
         & & \(K_2\) & 1.3944 & 1.3958 & 1.39514 & 1.3958 & 1.3958 & 1.3954 &  \\ 
        \cline{2-10}
         & \multirow{2}{*}{\(m_f\) [kg]} & \(K_1\) & 269.769 & 269.88 & 269.769 & 269.769 & 269.769 & 269.792 & \multirow{2}{*}{1.86853} \\ 
        \cline{3-9}
         & & \(K_2\) & 271.66 & 271.66 & 271.66 & 271.66 & 271.66 & 271.66 &  \\ 
        \hline
        \multirow{4}{*}{D} & \multirow{2}{*}{\(t_f\) [days]} & \(K_1\) & 25.756 & 25.744 & 25.746 & 25.735 & 25.754 & 25.747 & \multirow{2}{*}{-0.447361} \\ \cline{3-9}
         & & \(K_2\) & 25.527 & 24.706 & 25.422 & 24.817 & 26.026 & 25.300 &  \\ 
        \cline{2-10}
         & \multirow{2}{*}{\(m_f\) [kg]} & \(K_1\) & 946.697 & 946.697 & 946.697 & 946.696 & 946.697 & 946.697 & \multirow{2}{*}{0.029531} \\ 
        \cline{3-9}
         & & \(K_2\) & 946.701 & 946.719 & 946.761 & 946.7 & 946.75 & 946.726 &  \\ 
        \hline
        \multirow{4}{*}{E} & \multirow{2}{*}{\(t_f\) [days]} & \(K_1\) & 101.206 & 104.549 & 110.029 & 102.395 & 109.906 & 105.617 & \multirow{2}{*}{-10.5368} \\ \cline{3-9}
         & & \(K_2\) & 99.3076 & 91.6944 & 100.107 & 92.2472 & 92.0444 & 95.0801 &  \\ 
        \cline{2-10}
         & \multirow{2}{*}{\(m_f\) [kg]} & \(K_1\) & 1147.83 & 1145.12 & 1136.5 & 1130.6 & 1125.07 & 1137.02 & \multirow{2}{*}{89.5596} \\ 
        \cline{3-9}
         & & \(K_2\) & 1210.86 & 1238.72 & 1241.05 & 1235.96 & 1206.32 & 1226.58 &  \\ 
        \hline
    \end{tabularx}
\end{table}
\end{landscape}
\restoregeometry

\begin{figure}[h!]
    \centering
    \includegraphics[width=0.8\linewidth]{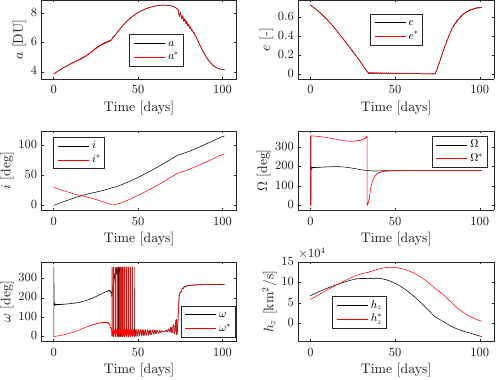}
    \caption{Case E orbital elements: minimum-time Q-law with \(\bm{K}_1\).}
    \label{fig: qlaw k1 oe}
\end{figure}

\begin{figure}[h!]
    \centering
    \includegraphics[width=0.8\linewidth]{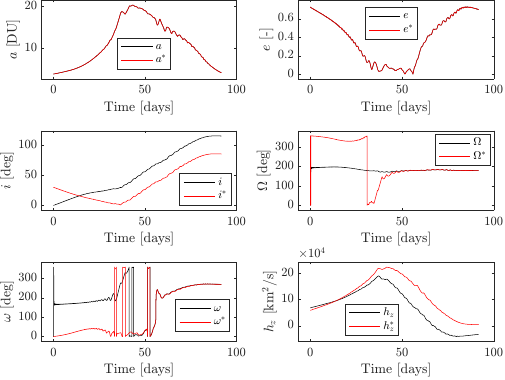}
    \caption{Case E orbital elements: minimum-time Q-law with \(\bm{K}_2\).}
    \label{fig: qlaw k2 oe}
\end{figure}

\begin{figure}[h!]
    \centering
    \includegraphics[width=0.8\linewidth]{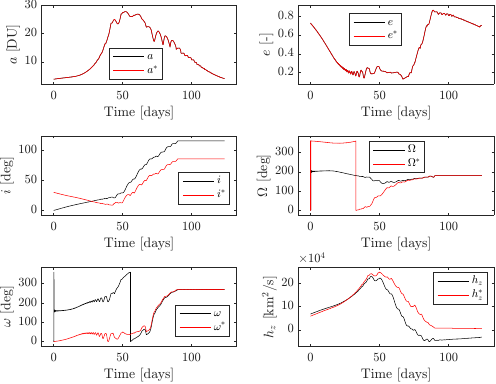}
    \caption{Case E orbital elements: minimum-fuel Q-law with \(\bm{K}_1\).}
    \label{fig: qlaw k1 mf oe}
\end{figure}

\begin{figure}[h!]
    \centering
    \includegraphics[width=0.8\linewidth]{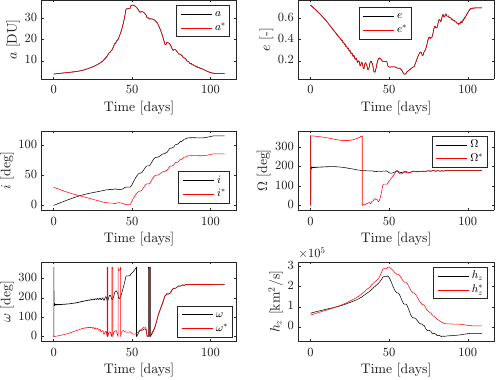}
    \caption{Case E orbital elements: minimum-fuel Q-law with \(\bm{K}_2\).}
    \label{fig: qlaw k2 mf oe}
\end{figure}

\begin{figure}[h!]
    \centering
    \begin{subfigure}[b]{0.49\textwidth}
        \centering
        \includegraphics[width=\textwidth]{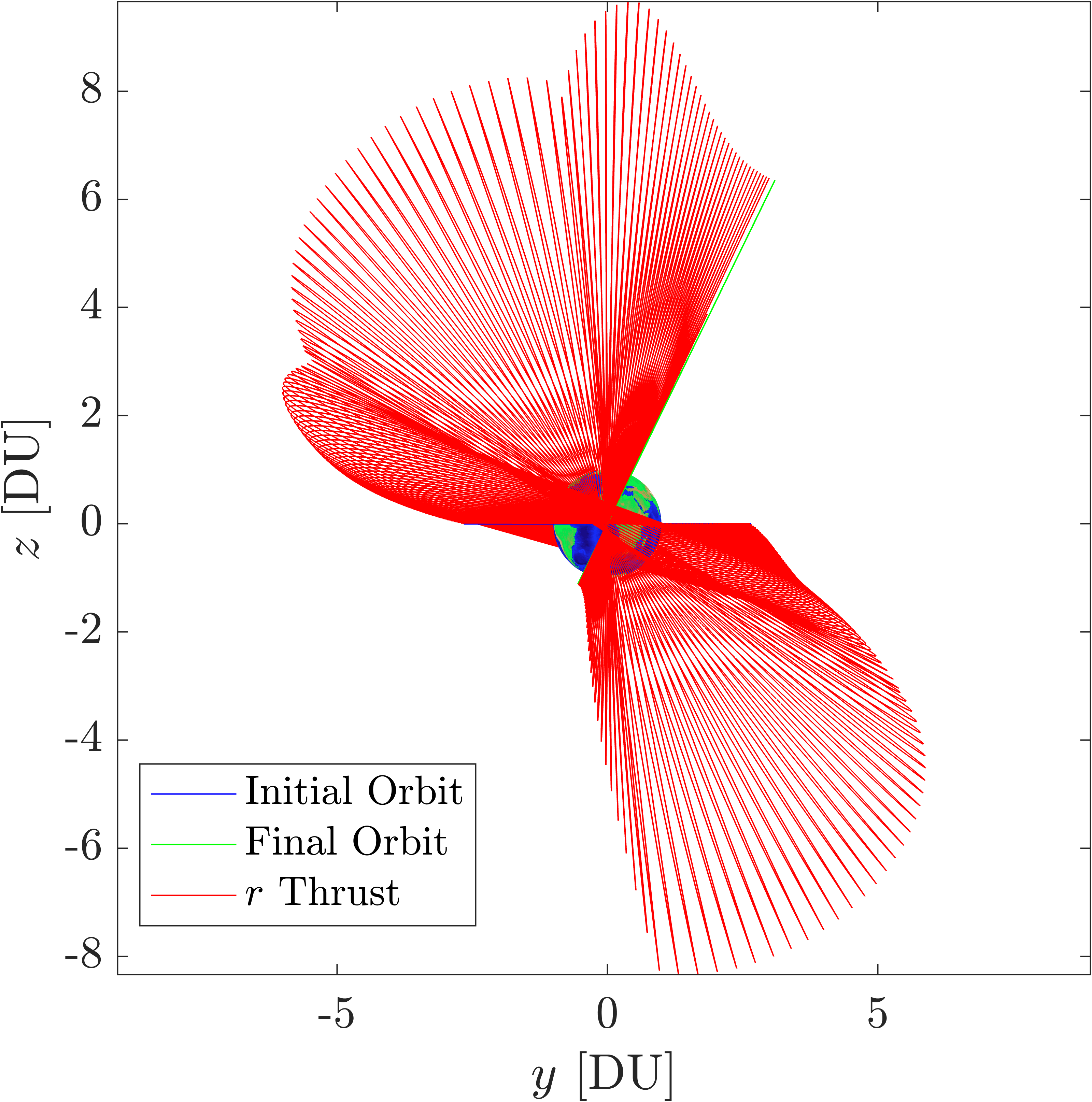}
        \caption{\(\bm{K}_1\).}
        \label{fig: qlaw k1 trajectory}
    \end{subfigure}
    \hfill
    \begin{subfigure}[b]{0.49\textwidth}
        \centering
        \includegraphics[width=\textwidth]{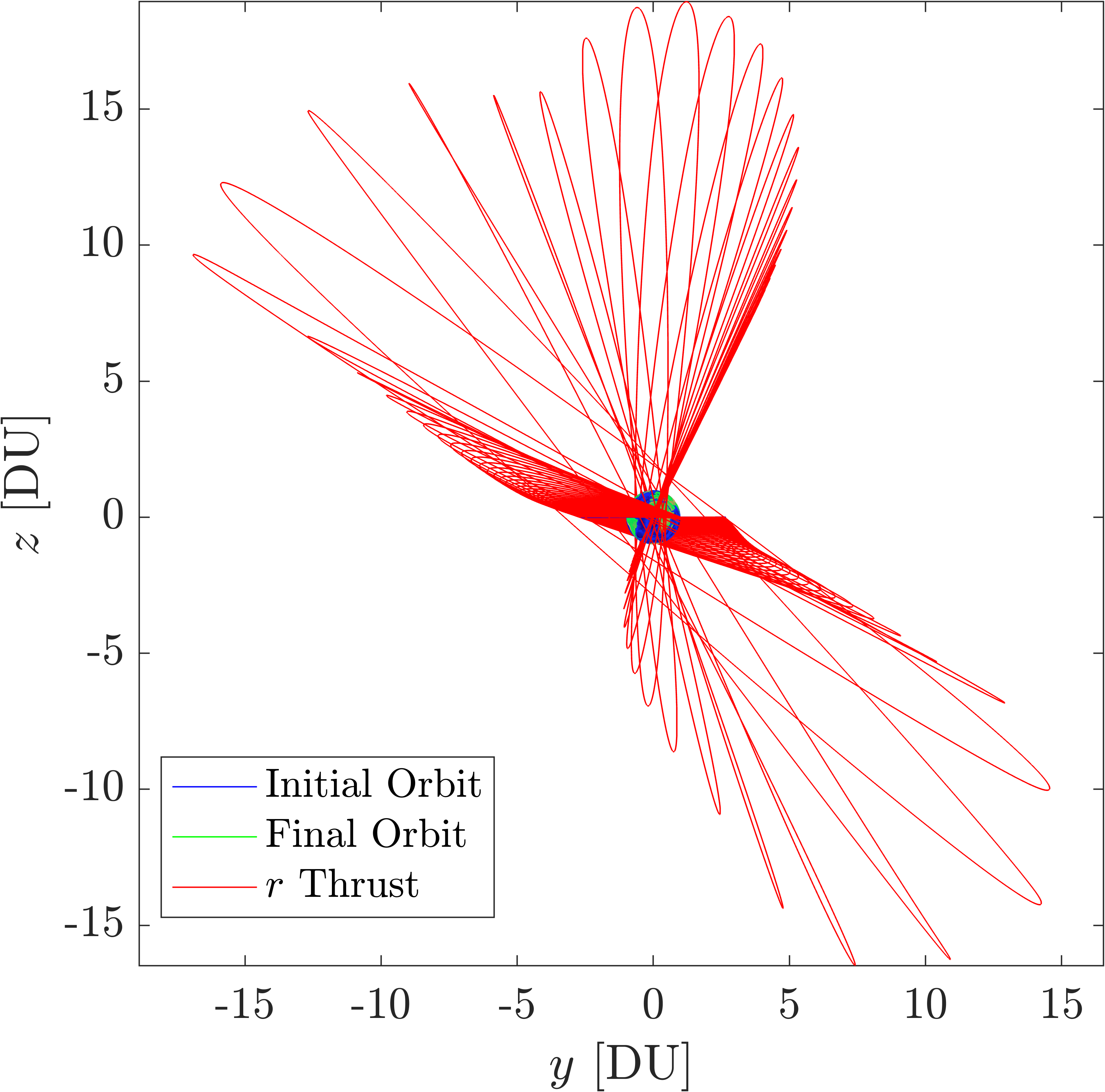}
        \caption{\(\bm{K}_2\)}
        \label{fig: qlaw k2 trajectory}
    \end{subfigure}
    \caption{Case E trajectory: minimum-time Q-law with \(\bm{K}_1\) and \(\bm{K}_2\).}
    \label{fig: qlaw case e mt trajectory}
\end{figure}

\begin{figure}[h!]
    \centering
    \begin{subfigure}[b]{0.49\textwidth}
        \centering
        \includegraphics[width=\textwidth]{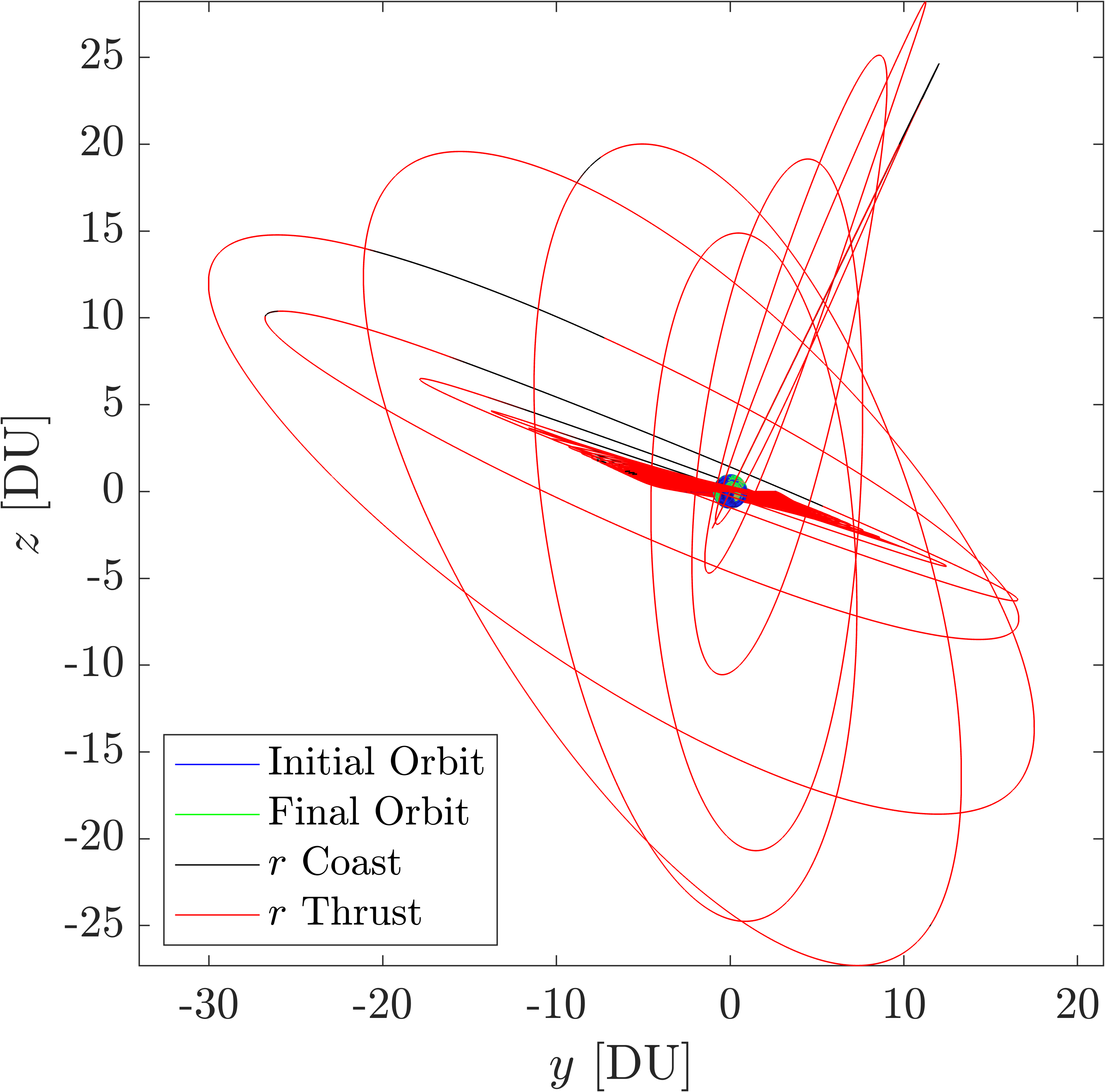}
        \caption{Trajectory projection onto \(z\)-\(y\) plane.}
        \label{fig: qlaw k1 mf trajectory normal view}
    \end{subfigure}
    \hfill
    \begin{subfigure}[b]{0.49\textwidth}
        \centering
        \includegraphics[width=\textwidth]{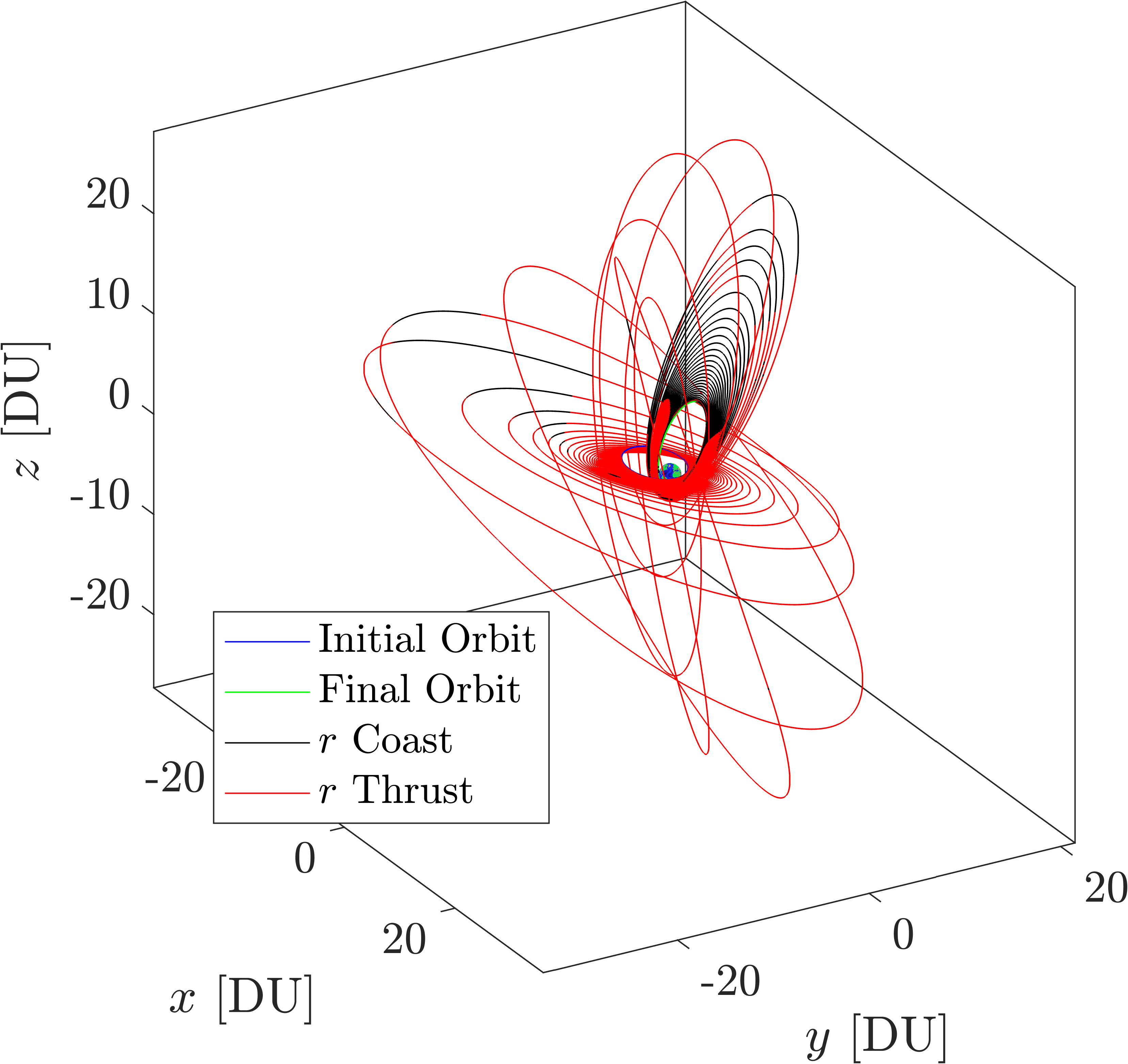}
        \caption{Three-dimensional view.}
        \label{qlaw k1 mf trajectory alt view}
    \end{subfigure}
    \caption{Case E trajectory: minimum-fuel Q-law with \(\bm{K}_1\).}
    \label{fig: qlaw k1 mf trajectory}
\end{figure}

\begin{figure}[h!]
    \centering
    \begin{subfigure}[b]{0.49\textwidth}
        \centering
        \includegraphics[width=\textwidth]{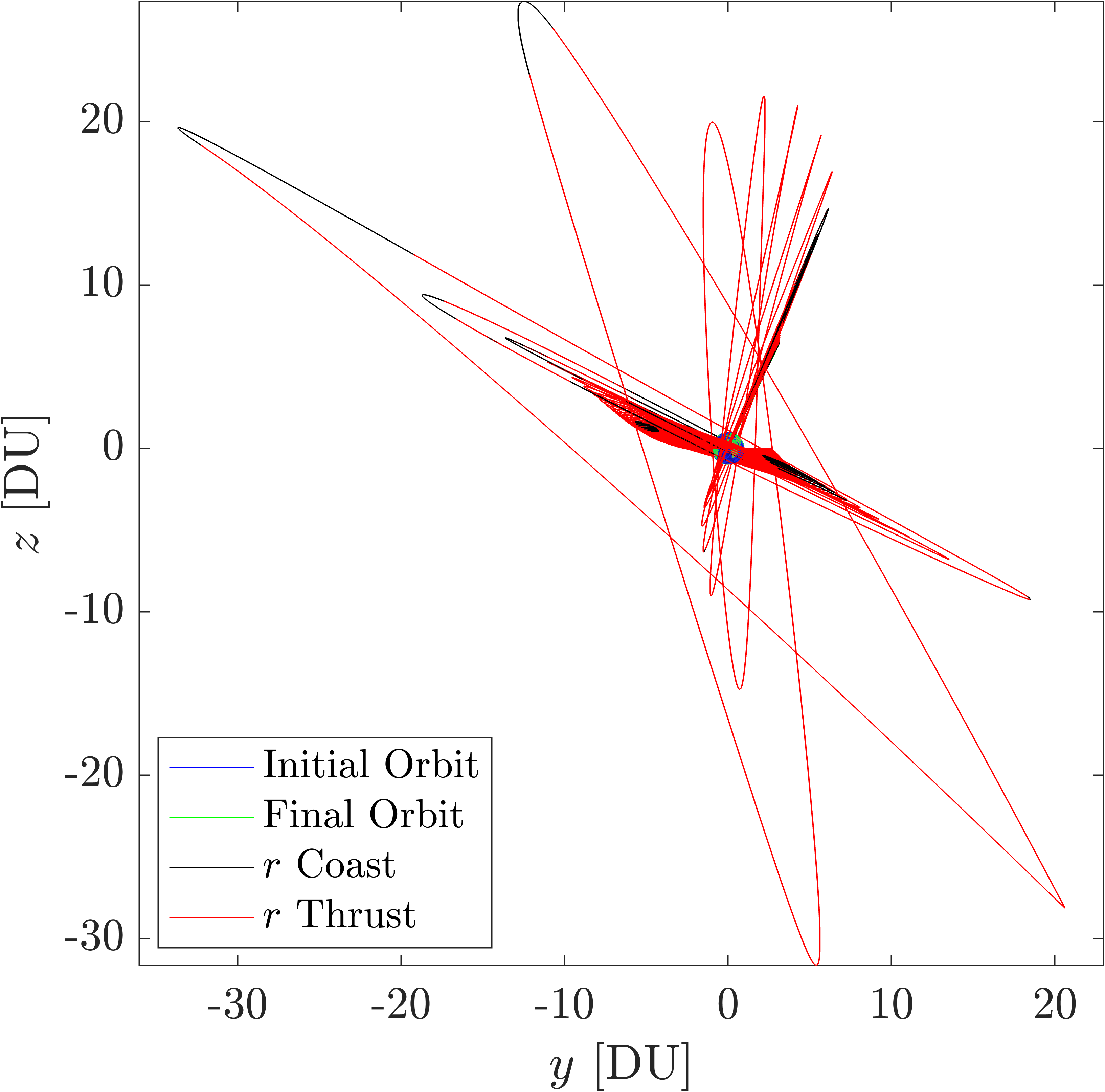}
        \caption{Trajectory projection onto \(z\)-\(y\) plane.}
        \label{fig: qlaw k2 mf trajectory normal view}
    \end{subfigure}
    \hfill
    \begin{subfigure}[b]{0.49\textwidth}
        \centering
        \includegraphics[width=\textwidth]{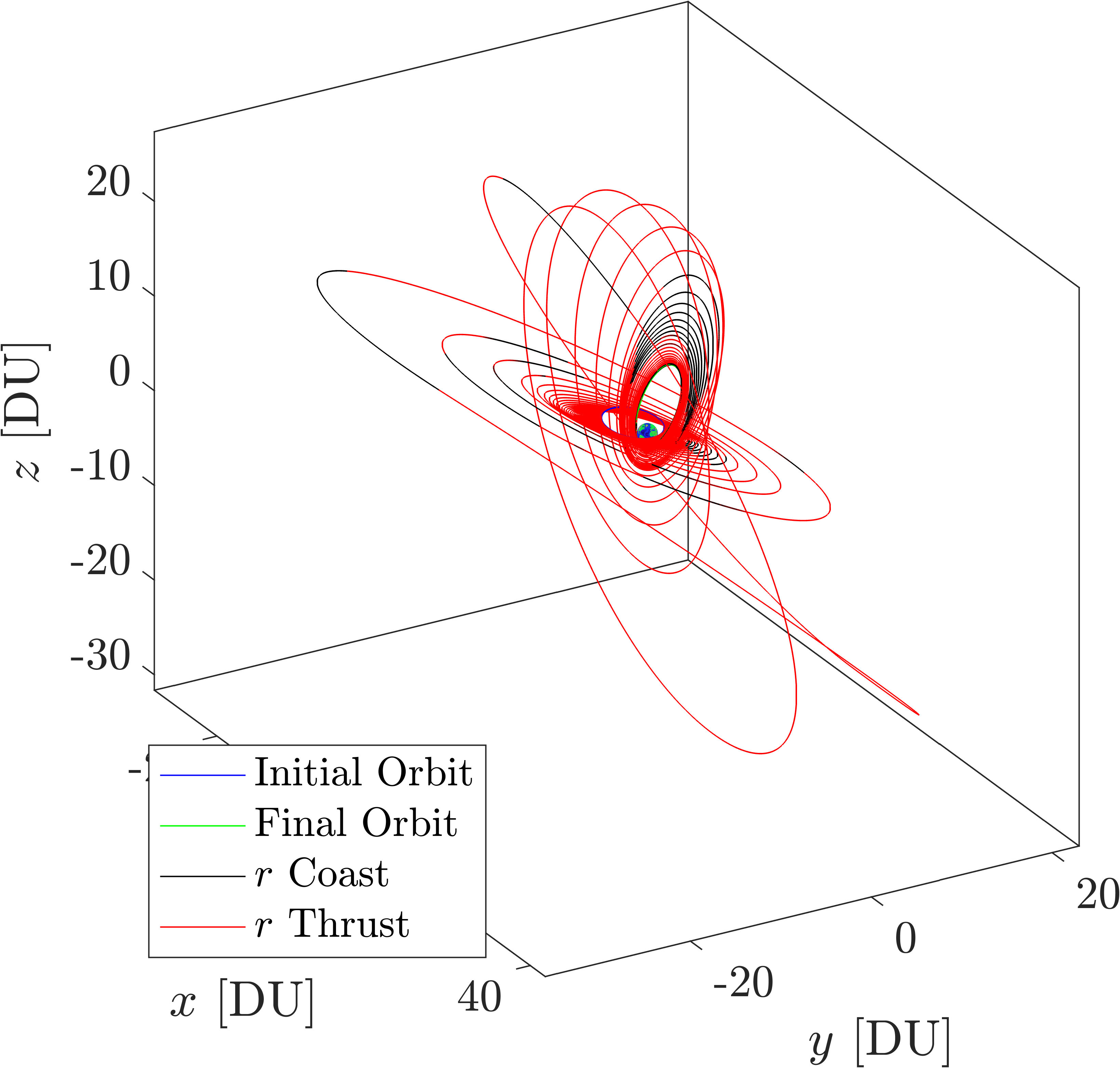}
        \caption{Three-dimensional view.}
        \label{qlaw k2 mf trajectory alt view}
    \end{subfigure}
    \caption{Case E trajectory: minimum-fuel Q-law with \(\bm{K}_2\).}
    \label{fig: qlaw k2 mf trajectory}
\end{figure}

Figures \ref{fig: qlaw k1 oe}, \ref{fig: qlaw k2 oe}, \ref{fig: qlaw k1 mf oe}, and \ref{fig: qlaw k2 mf oe} show the orbital elements time histories for the best Case E minimum-time and minimum-fuel Q-law solutions for each penalty matrix type. Figures \ref{fig: qlaw case e mt trajectory}, \ref{fig: qlaw k1 mf trajectory}, and \ref{fig: qlaw k2 mf trajectory} show the trajectories for these solutions. Note that the rotated Case E* boundary conditions are plotted in red in the orbital element plots along with the original Case E orbital elements in black. The problem is solved in the Case E* frame and the solution is simply rotated back into the Case E frame with single rotation about the \(x\)-axis. The trajectories show the Case E trajectory solution. Similar to the results for the generic Lyapunov control law in the previous section, the full penalty matrix \(\bm{K}_2\) changes orbital elements over time quite differently than with the diagonal penalty matrix, \(\bm{K}_1\). Note that any oscillations or sharp jumps in \(\Omega\) and \(\omega\) between 0\(^\circ\) and 360\(^\circ\) are due to the orbit being nearly circular and/or equatorial leading to ambiguity in \(\Omega\) and \(\omega\).

\section{Conclusion} \label{sec: conclusion}
Parameterization and optimization of positive-definite penalty matrices has widespread applications in modern control and optimization theories. Diagonal parameterization is a common approach, but it represents only a subset of the full solution space, as it neglects cross-coupled terms. We propose an efficient parameterization that is based on eigendecomposition of symmetric matrices, in which the set of parameters consists of eigenvalues and variables used to parameterize the orthogonal eigenvector matrices. Considering a full penalty matrix can broaden the solution space by including cross-coupling error terms in quadratic costs that are used extensively for designing modern controls, such as a linear quadratic regulator (LQR) or Lyapunov-based nonlinear controllers, as we have demonstrated through the results.

Several control problems were solved with Lyapunov-based control laws and LQR for the standard diagonal and the proposed full penalty matrix representations. We clarify that solutions are by no means optimal due to limited runs of PSO and limited computational resources. Additionally, other global optimization algorithms may even perform better than PSO. For instance, Ref. \cite{lee_design_2005} finds more optimal transfer solutions using Q-law by performing extensive parameter searches with genetic algorithm over multiple days of CPU time. The matter of optimizing the parameters should be investigated further. However, our comparative results demonstrate that improved solutions can be obtained with the full penalty matrix with respect to the considered optimality criteria of each problem. The improvement is especially evident for classes of more nonlinear and long-duration problems for which optimization of the off-diagonal terms broadens the optimality of the solutions (e.g., large changes in orbital elements for the low-thrust trajectories). In the most extreme case, the improvement in performance was 65\%. Our hope is that the simplicity of the demonstrated methods allows researchers/practitioners in many fields of control theory and optimization to potentially obtain more optimal solutions to their practical control problems.


\appendix
\section{Numerical Values of the Penalty Matrices}
\label{app1}

For reproducibility, we provide the penalty matrices for some of the solutions we obtained in this work. Eqs.~\eqref{eq: best of lqr detumbling} and \eqref{eq: best of detumbling lc} show the penalty matrices that were found for the best detumbling attitude control solution for LQR and Lyapunov-based control, respectively, corresponding with the solutions shown in Figures \ref{fig: lqr manuever 1 plots} and \ref{fig: lyapunov based manuever 1 plots}. Eqs.~\eqref{eq: glc low thrust case E}, \eqref{eq: best of qlaw min time}, \eqref{eq: best of qlaw min fuel} show the penalty matrices for the Case E for the minimum-time generic LC solution, minimum-time Q-law solution, and minimum-fuel Q-law solution, respectively.

\begin{subequations} \label{eq: best of lqr detumbling}
    \small
    \begin{align}
        \bm{Q}_1 & = \text{diag}\left(1.05309,0.00299126,0.68307,4.5293,9.87489,0.0452255\right) \\
        \bm{R}_1 & = \text{diag}\left(8.36281,4.40636,9.23266\right) \\
        \bm{Q}_2 & = \left[\begin{matrix} 1.185 & -0.613646 & 0.311778 & 0.48939 & 0.343032 & 0.0979317 \\ 
            -0.613646 & 0.698828 & 0.572306 & -1.82684 & -0.0230107 & -0.0140387 \\ 
            0.311778 & 0.572306 & 1.84127 & -3.45008 & 0.448065 & 0.110059 \\ 
            0.48939 & -1.82684 & -3.45008 & 7.57063 & -0.592071 & -0.132696 \\ 
            0.343032 & -0.0230107 & 0.448065 & -0.592071 & 0.172537 & 0.0455803 \\ 
            0.0979317 & -0.0140387 & 0.110059 & -0.132696 & 0.0455803 & 0.0121875 \end{matrix}\right] \\
        \bm{R}_2 & = \begin{bmatrix} 
            9.57947 & 1.63521 & 0.963686 \\ 
            1.63521 & 3.39122 & -4.30809 \\ 
            0.963686 & -4.30809 & 6.53455
            \end{bmatrix}
    \end{align}
\end{subequations}

\begin{subequations}\label{eq: best of detumbling lc}
    \begin{align}
        \bm{K}_{p,1} & = \text{diag}\left(0.192456,0.167666,0.272458\right) \\ 
        \bm{K}_{d,1} & = \text{diag}\left(1.82655,2.25069,3.04384\right) \\
        \bm{K}_{p,2} & = \begin{bmatrix} 
        2.2675 & -0.295543 & 1.97134 \\ 
        -0.295543 & 0.915225 & -0.761434 \\ 
        1.97134 & -0.761434 & 2.35022
        \end{bmatrix} \\
        \bm{K}_{d,2} & = \begin{bmatrix} 
        7.03549 & 0.367325 & 1.74591 \\ 
        0.367325 & 7.73374 & -2.65179 \\ 
        1.74591 & -2.65179 & 2.34927
        \end{bmatrix}
    \end{align}
\end{subequations}

\begin{subequations}\label{eq: glc low thrust case E}
    \begin{align}
        \bm{K}_1 & = \text{diag}\left(6.5225,98.1494,8.9658,10.2037,97.3815,20.7142\right) \\
        \bm{K}_2 & = \left[\begin{matrix} 39.4746 & 17.5941 & -2.0538 & -3.3242 & 0.1723 & -0.3125 \\
                                         17.5941 & 68.1822 & -3.7857 & -6.9744 & 0.1840 & -0.5589 \\
                                         -2.0538 & -3.7857 & 4.6930 & 0.5193 & 2.9905 & 0.2195 \\
                                         -3.3242 & -6.9744 & 0.5193 & 11.5512 & 0.4138 & -5.6421 \\
                                         0.1723 & 0.1840 & 2.9905 & 0.4138 & 77.1079 & 1.5671 \\
                                         -0.3125 & -0.5589 & 0.2195 & -5.6421 & 1.5671 & 81.3064 \end{matrix}\right]
    \end{align}
\end{subequations}

\begin{subequations}\label{eq: best of qlaw min time}
    \begin{align}
        \bm{K}_1 & = \text{diag}\left(2.54742,0.00530434,0.242459,9.64791,7.76675 \right) \\
        \bm{K}_2 & = \left[\begin{matrix} 9.61437 & 0.59816 & 0.727462 & 0.0886329 & 0.0288422 \\ 
                                            0.59816 & 5.65613 & -4.0757 & -1.62804 & -1.24825 \\ 
                                            0.727462 & -4.0757 & 3.52475 & 0.90911 & 1.25853 \\ 
                                            0.0886329 & -1.62804 & 0.90911 & 4.76366 & 0.652616 \\ 
                                            0.0288422 & -1.24825 & 1.25853 & 0.652616 & 2.68927 \end{matrix}\right]
    \end{align}
\end{subequations}

\begin{subequations}\label{eq: best of qlaw min fuel}
    \begin{align}
        \bm{K}_1 & = \text{diag}\left(8.65871,2.22226,4.21207,4.84089,8.3779\right) \\
        \bm{K}_2 & = \left[\begin{matrix} 8.25151 & 3.14031 & 1.22609 & 0.337485 & 0.0941679\\ 
                                        3.14031 & 3.62918 & -1.34696 & 0.0247984 & 0.066369 \\ 
                                        1.22609 & -1.34696 & 4.13089 & 0.256997 & -0.0302524 \\ 
                                        0.337485 & 0.0247984 & 0.256997 & 8.95998 & -0.429254 \\ 
                                        0.0941679 & 0.066369 & -0.0302524 & -0.429254 & 7.19688 \end{matrix}\right]
    \end{align}
\end{subequations}

\bibliographystyle{elsarticle-num} 
\bibliography{references}







\end{document}